# STATISTICAL EIGEN-INFERENCE FROM LARGE WISHART MATRICES

By N. Raj Rao,[1,2] James A. Mingo,[3] Roland Speicher[3,4]
and Alan Edelman[1]

*Massachusetts Institute of Technology, Queen's University,
Queen's University and Massachusetts Institute of Technology*

We consider settings where the observations are drawn from a zero-mean multivariate (real or complex) normal distribution with the population covariance matrix having eigenvalues of arbitrary multiplicity. We assume that the eigenvectors of the population covariance matrix are unknown and focus on inferential procedures that are based on the sample eigenvalues alone (i.e., "eigen-inference").

Results found in the literature establish the asymptotic normality of the fluctuation in the trace of powers of the sample covariance matrix. We develop concrete algorithms for analytically computing the limiting quantities and the covariance of the fluctuations. We exploit the asymptotic normality of the trace of powers of the sample covariance matrix to develop eigenvalue-based procedures for testing and estimation. Specifically, we formulate a simple test of hypotheses for the population eigenvalues and a technique for estimating the population eigenvalues in settings where the cumulative distribution function of the (nonrandom) population eigenvalues has a staircase structure.

Monte Carlo simulations are used to demonstrate the superiority of the proposed methodologies over classical techniques and the robustness of the proposed techniques in high-dimensional, (relatively) small sample size settings. The improved performance results from the fact that the proposed inference procedures are "global" (in a sense that we describe) and exploit "global" information thereby overcoming the inherent biases that cripple classical inference procedures which are "local" and rely on "local" information.

Received January 2007; revised December 2007.

[1]Supported in part by NSF Grant DMS-04-11962.

[2]Supported in part by an Office of Naval Research Special Post-Doctoral Award Grant N00014-07-1-0269.

[3]Supported in part by Discovery Grants and a Leadership Support Initiative Award from the Natural Sciences and Engineering Research Council of Canada.

[4]Supported in part by a Premier's Research Excellence Award from the Province of Ontario and a Killam Fellowship from the Canada Council for the Arts.

*AMS 2000 subject classifications.* 62510, 62E20, 15A52.

*Key words and phrases.* Sample covariance matrices, random matrix theory, Wishart matrices, second order freeness, free probability, eigen-inference, linear statistics.

**1. Introduction.** Let $\mathbf{X} = [\mathbf{x}_1, \ldots, \mathbf{x}_n]$ be a $p \times n$ data matrix where $\mathbf{x}_1, \ldots, \mathbf{x}_n$ denote $n$ independent measurements, where for each $i$, $\mathbf{x}_i$ has a $p$-dimensional (real or complex) Gaussian distribution with mean zero, and positive definite covariance matrix $\boldsymbol{\Sigma}$. When the samples are complex, the real and imaginary components are assumed to be independent, identically distributed zero-mean Gaussian vectors with a covariance of $\boldsymbol{\Sigma}/2$. The sample covariance matrix (SCM) when formed from these $n$ samples as

$$(1.1) \qquad \mathbf{S} := \frac{1}{n} \sum_{i=1}^{n} \mathbf{x}_i \mathbf{x}_i' = \frac{1}{n} \mathbf{X} \mathbf{X}',$$

with $'$ denoting the conjugate transpose, has the (central) Wishart distribution [Wishart (1928)]. We focus on inference problems for parameterized covariance matrices modeled as $\boldsymbol{\Sigma_\theta} = \mathbf{U} \boldsymbol{\Lambda_\theta} \mathbf{U}'$ where

$$(1.2) \qquad \boldsymbol{\Lambda_\theta} = \begin{bmatrix} a_1 \mathbf{I}_{p_1} & & & \\ & a_2 \mathbf{I}_{p_2} & & \\ & & \ddots & \\ & & & a_k \mathbf{I}_{p_k} \end{bmatrix},$$

where $a_1 > \cdots > a_k$ and $\sum_{j=1}^{k} p_j = p$. Defining $t_i = p_i/p$ allows us to conveniently express the $2k-1$-dimensional parameter vector as $\boldsymbol{\theta} = (t_1, \ldots, t_{k-1}, a_1, \ldots, a_k)$ with the obvious nonnegativity constraints on the elements.

Models of the form in (1.2) arise as a special case whenever the measurements are of the form

$$(1.3) \qquad \mathbf{x}_i = \mathbf{A} \mathbf{s}_i + \mathbf{z}_i \qquad \text{for } i = 1, \ldots, n,$$

where $\mathbf{z}_i \sim \mathcal{N}_p(0, \boldsymbol{\Sigma}_z)$ denotes a $p$-dimensional (real or complex) Gaussian noise vector with covariance $\boldsymbol{\Sigma}_z$, $\mathbf{s}_i \sim \mathcal{N}_k(\mathbf{0}, \boldsymbol{\Sigma}_s)$ denotes a $k$-dimensional zero-mean (real or complex) Gaussian signal vector with covariance $\boldsymbol{\Sigma}_s$, and $\mathbf{A}$ is a $p \times k$ unknown nonrandom matrix. In array processing applications, the $j$th column of the matrix $\mathbf{A}$ encodes the parameter vector associated with the $j$th signal whose amplitude is described by the $j$th element of $\mathbf{s}_i$. See, for example, the text by Van Trees (2002).

Since the signal and noise vectors are independent of each other, the covariance matrix of $\mathbf{x}_i$ can be decomposed as

$$(1.4) \qquad \boldsymbol{\Sigma} = \boldsymbol{\Psi} + \boldsymbol{\Sigma}_z,$$

where $\boldsymbol{\Sigma}_z$ is the covariance of $\mathbf{z}$ and $\boldsymbol{\Psi} = \mathbf{A} \boldsymbol{\Sigma}_s \mathbf{A}'$. One way of obtaining $\boldsymbol{\Sigma}$ with eigenvalues of the form in (1.2) is when $\boldsymbol{\Sigma}_z = \sigma^2 \mathbf{I}$ so that the $n - k$ smallest eigenvalues of $\boldsymbol{\Sigma}$ are equal to $\sigma^2$. Then, if the matrix $\mathbf{A}$ is of full column rank and the covariance matrix of the signals $\boldsymbol{\Sigma}_s$ is nonsingular, the $p - k$ (with $k < p$ here) smallest eigenvalues of $\boldsymbol{\Psi}$ are equal to zero so that the eigenvalues of $\boldsymbol{\Sigma}$ will be of the form in (1.2). Alternatively, if the eigenvalues of $\boldsymbol{\Psi}$ and $\boldsymbol{\Sigma}_z$ have the identical subspace structure, that is, in (1.2), $t_i^\Psi = t_i^{\Sigma_z}$ for all $i$, then whenever the eigenvectors associated with each of the subspaces of $\boldsymbol{\Psi}$ and $\boldsymbol{\Sigma}_z$ align, the eigenvalues of $\boldsymbol{\Sigma}$ will have the subspace structure in (1.2).

Additionally, from an identifiability point of view, as shall be discussed in Section 7, if the practitioner has reason to believe that the population eigenvalues are organized in several clusters about $a_i \pm a_i \sqrt{p/n}$, then the use of the model in (1.2) with a block subspace structure will also be justified.

1.1. *Inferring the population eigenvalues from the sample eigenvalues.* While inference problems for these models have been documented in texts such as (Muirhead (1982)), the inadequacies of classical algorithms in high-dimensional, (relatively) small sample size settings have not been adequately addressed. We highlight some of the prevalent issues in the context of statistical inference and hypothesis testing.

In the landmark paper [Anderson (1963)], the theory was developed that describes the (large sample) asymptotics of the sample eigenvalues (in the real-valued case) for such models when the true covariance matrix has eigenvalues of arbitrary multiplicity. Indeed, for arbitrary covariance $\boldsymbol{\Sigma}$, the joint density function of the eigenvalues $l_1, \ldots, l_p$ of the SCM $\mathbf{S}$ when $n > p + 1$ is shown to be given by

$$(1.5) \quad \widetilde{Z}_{p,n}^\beta \sum_{i=1}^p l_i^{\beta(n-p+1)/2-1} \prod_{i<j}^p |l_i - l_j|^\beta \int_{\mathbf{Q}} \exp\left(-\frac{n\beta}{2} \operatorname{Tr}(\boldsymbol{\Sigma}^{-1} \mathbf{V} \mathbf{S} \mathbf{V}')\right) d\mathbf{V},$$

where $l_1 > \cdots > l_p > 0$, $\widetilde{Z}_{p,n}^\beta$ is a normalization constant, and $\beta = 1$ (or 2) when $\mathbf{S}$ is real (resp., complex). In (1.5), $\mathbf{Q} \in \mathbf{O}(p)$ when $\beta = 1$ while $\mathbf{Q} \in \mathbf{U}(p)$ when $\beta = 2$ where $\mathbf{O}(p)$ and $\mathbf{U}(p)$ are, respectively, the set of $p \times p$ orthogonal and unitary matrices with Haar measure. Anderson notes that

> If the characteristic roots of $\boldsymbol{\Sigma}$ are different, the deviations ... from the corresponding population quantities are asymptotically normally distributed. When some of the roots of $\boldsymbol{\Sigma}$ are equal, the asymptotic distribution cannot be described so simply.

Indeed, the difficulty alluded to $\boldsymbol{\Sigma}$ arises due to the presence of the integral over orthogonal (or unitary) group on the right-hand side of (1.5). This problem is compounded in situations when some of the eigenvalues of $\Sigma$ are equal as is the case for the model considered in (1.2). In such settings, large sample approximations for this multivariate integral have been used [see, e.g., Muirhead (1982), page 403, Corollary 9.5.6, Butler and Wood (2002, 2005)]. For the problem of interest, Anderson uses just such an approximation to derive the maximum-likelihood estimate of the population eigenvalues, $a_l$, as

$$(1.6) \qquad \widehat{a}_l \approx \frac{1}{p_l} \sum_{j \in N_l} \widehat{\lambda}_j \qquad \text{for } l = 1, \dots, k,$$

where $\widehat{\lambda}_j$ are the sample eigenvalues (arranged in descending order) and $N_l$ is the set of integers $p_1 + \cdots + p_{l-1} + 1, \dots, p_1 + \cdots + p_l$. This is a reasonable estimator that works well in practice when $n \gg p$. The large sample size asymptotics are, however, of limited utility because they ignore the (significant) effect of the dimensionality of the system on the behavior of the sample eigenvalues.

Consequently, (large sample size) asymptotic predictions, derived under the $p$ fixed, $n \to \infty$ regime do not account for the additional complexities that arise in situations where the sample size $n$ is large but the dimensionality $p$ is of comparable order. Furthermore, the estimators developed using the classical large sample asymptotics invariably become degenerate whenever $p > n$, so that $p - n$ of the sample eigenvalues will identically equal to zero. For example, when $n = p/2$, and there are two distinct population eigenvalues each with multiplicity $p/2$, then the estimate of the smallest eigenvalue using (1.6) will be zero. Other such scenarios where the population eigenvalue estimates obtained using (1.6) are meaningless are easy to construct and are practically relevant in many applications such as radar and sonar signal processing, and many more. See, for example, the text by Van Trees (2002) and the work of Smith (2005).

There are, of course, other strategies one may employ for inferring the population eigenvalues. One might consider a maximum-likelihood technique based on maximizing the log-likelihood function of the observed data $\mathbf{X}$ which is given by (ignoring constants)

$$l(\mathbf{X}|\mathbf{\Sigma}) := -n(\text{Tr}\,\mathbf{S}\mathbf{\Sigma}^{-1} + \log \det \mathbf{\Sigma}),$$

or, equivalently, when $\mathbf{\Sigma} = \mathbf{U}\mathbf{\Lambda}\mathbf{U}'$, by minimizing the objective function

$$(1.7) \qquad h(\mathbf{X}|\mathbf{U}, \mathbf{\Lambda}) = (\text{Tr}\,\mathbf{S}\mathbf{U}\mathbf{\Lambda}^{-1}\mathbf{U}' + \log \det \mathbf{\Lambda}).$$

What should be apparent on inspecting (1.7) is that the maximum-likelihood estimation of the parameters of $\mathbf{\Lambda}$ of the form in (1.2) requires us to model the population eigenvectors $\mathbf{U}$ as well (except when $k = 1$). If $\mathbf{U}$ were known a priori, then an estimate of $a_l$ obtained as

$$(1.8) \qquad \widehat{a}_l \approx \frac{1}{p_l} \sum_{j \in N_l} (\mathbf{U}'\mathbf{S}\mathbf{U})_{j,j} \qquad \text{for } l = 1, \dots, k,$$

where $N_l$ is the set of integers $p_1 + \cdots + p_{l-1} + 1, \dots, p_1 + \cdots + p_l$, will provide a good estimate. In practical applications, the population eigenvectors might either be unknown or be misspecified leading to faulty inference. Hence it is important to have the ability to perform statistically sound, computationally feasible *eigen-inference* of the population eigenvalues, that is, from the sample eigenvalues alone, in a manner that is robust to high-dimensionality and sample size constraints.

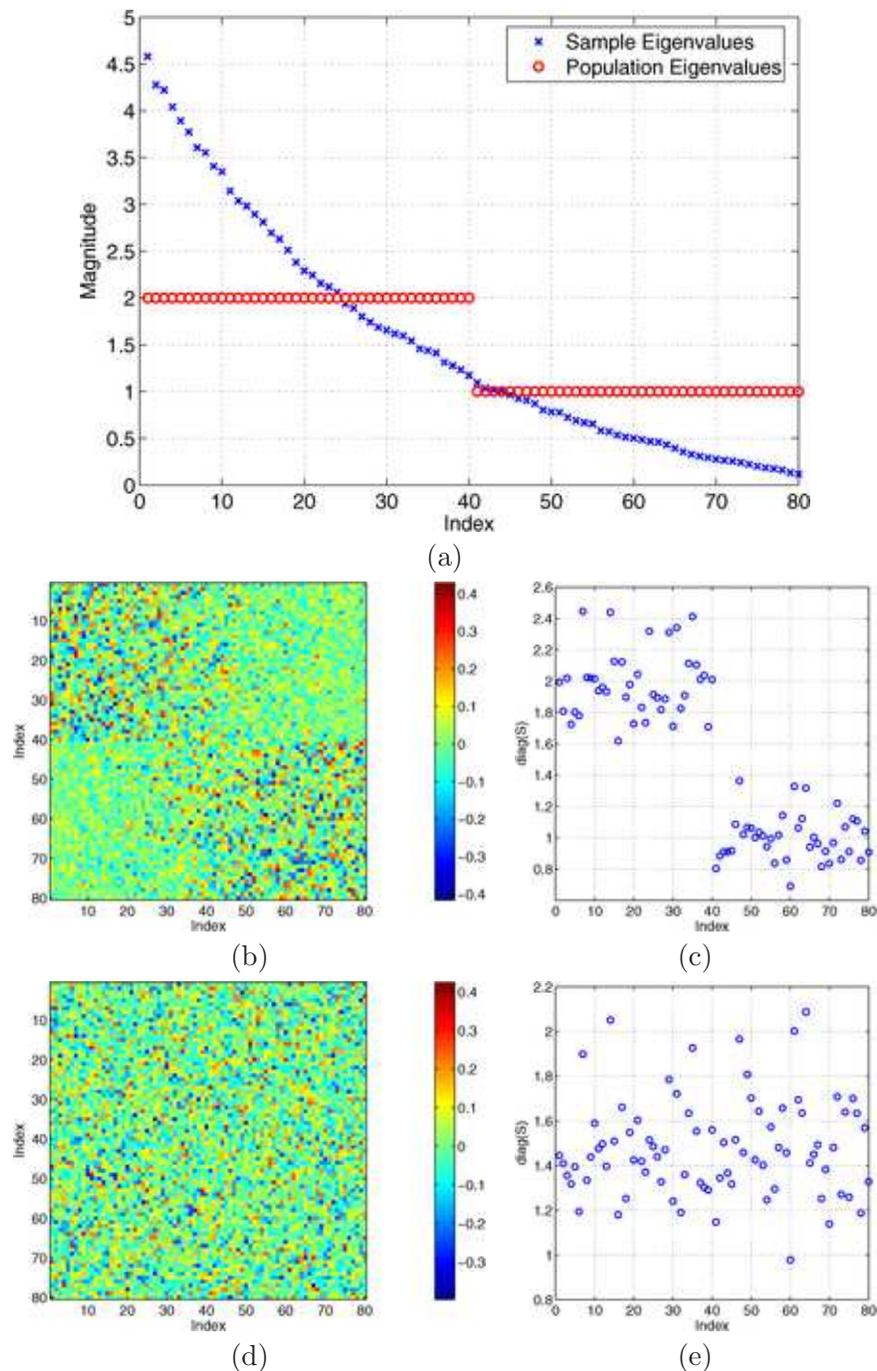

FIG. 1. *The challenge of estimating the population eigenvalues from the sample eigenvalues in high-dimensional settings.* (a) *Sample eigenvalues versus true eigenvalues* ($p = 80$, $n = 160$). (b) *Sample eigenvectors when* $\mathbf{U} = \mathbf{I}$. (c) *Diagonal elements of* $\mathbf{S}$ *when* $\mathbf{U} = \mathbf{I}$. (d) *Sample eigenvectors for arbitrary* $\mathbf{U}$. (e) *Diagonal elements of* $\mathbf{S}$ *for arbitrary* $\mathbf{U}$.

We illustrate the difficulties encountered in high-dimensional settings with an example (summarized in Figure 1) of a SCM constructed from a covariance matrix modeled as $\boldsymbol{\Sigma} = \mathbf{U}\boldsymbol{\Lambda}\mathbf{U}'$ with $p = 80$ and sample size $n = 160$. Half of the eigenvalues of $\boldsymbol{\Lambda}$ are equal to 2 while the remainder are equal to 1. The sample eigenvalues are significantly *blurred*, relative to the true eigenvalues as shown in Figure 1(a). Figure 1(b) and (d) plot the sample eigenvectors for the case when the true eigenvectors $\mathbf{U} = \mathbf{I}$, and an arbitrary $\mathbf{U}$, respectively. Figure 1(c) and (e) plot the diagonal elements $(\mathbf{S})_{j,j}$. Thus, if the true eigenvector was indeed $\mathbf{U} = \mathbf{I}$, then an estimate of the population eigenvalues formed as in (1.8) yields a good estimate; when $\mathbf{U} \neq \mathbf{I}$, however, the estimate is very poor.

1.2. *Testing for equality of population eigenvalues.* Similar difficulties are encountered in problems of testing as well. In such situations, for testing the hypothesis

$$\lambda_{p_1 + \cdots + p_{l-1} + 1} = \lambda_{p_1 + \cdots + p_{l-1} + 1, \ldots, p_1 + \cdots + p_l}.$$

Anderson proposes the likelihood criterion given by

$$(1.9) \qquad V_l = \left[ \prod_{j \in N_l} \widehat{\lambda}_j \Big/ \left( p_k^{-1} \sum_{j \in N_l} \widehat{\lambda}_j \right)^{p_k} \right]^{n/2} \qquad \text{for } l = 1, \ldots, k,$$

where $\widehat{\lambda}_j$ are the sample eigenvalues (arranged in descending order) and, as before, $N_l$ is the set of integers $p_1 + \cdots + p_{l-1} + 1, \ldots, p_1 + \cdots + p_l$. The test in (1.9) suffers from the same deficiency as the population eigenvalue estimator in (1.6)—it becomes degenerate when $p > n$. When the population eigenvectors $\mathbf{U}$ are known, (1.9) may be modified by forming the criterion

$$(1.10) \qquad \left[ \prod_{j \in N_l} (\mathbf{U}'\mathbf{S}\mathbf{U})_{j,j} \Big/ \left( p_k^{-1} \sum_{j \in N_l} (\mathbf{U}'\mathbf{S}\mathbf{U})_{j,j} \right)^{p_k} \right]^{n/2} \qquad \text{for } l = 1, \ldots, k.$$

When the eigenvectors are misspecified, the inference provided will be faulty. For the earlier example, Figure 1(e) illustrates this for the case when it is assumed that the population eigenvectors are $\mathbf{I}$ when they are really $\mathbf{U} \neq \mathbf{I}$. Testing the hypothesis $\boldsymbol{\Sigma} = \boldsymbol{\Sigma}_0$ reduces to testing the hypothesis $\boldsymbol{\Sigma} = \mathbf{I}$, given samples $\widetilde{\mathbf{x}}_i$ for $i = 1, \ldots, n$, where $\widetilde{\mathbf{x}}_i = \boldsymbol{\Sigma}_0^{-1/2}\mathbf{x}_i$. The robustness of tests for sphericity in high-dimensional settings has been extensively discussed by Ledoit and Wolf (2002) and is the focus of some recent work by Srivastava (2005, 2006).

TABLE 1
*Structure of proposed algorithms*

| Testing: | $H_{\boldsymbol{\theta}_0} : h(\boldsymbol{\theta}) := \mathbf{v}_{\boldsymbol{\theta}}^T \mathbf{Q}_{\boldsymbol{\theta}}^{-1} \mathbf{v}_{\boldsymbol{\theta}},$ | Recommend $\dim(\mathbf{v}_{\boldsymbol{\theta}}) = 2$ |
|---|---|---|
| Estimation: | $\widehat{\boldsymbol{\theta}} = \underset{\boldsymbol{\theta} \in \Theta}{\arg\min} \{ \mathbf{v}_{\boldsymbol{\theta}}^T \mathbf{Q}_{\boldsymbol{\theta}}^{-1} \mathbf{v}_{\boldsymbol{\theta}} + \log \det \mathbf{Q}_{\boldsymbol{\theta}} \},$ | Recommend $\dim(\mathbf{v}_{\boldsymbol{\theta}}) = \dim(\boldsymbol{\theta}) + 1$ |
| Legend: | $(\mathbf{v}_{\boldsymbol{\theta}})_j = p \times \left( \frac{1}{p} \operatorname{Tr} S^j - \alpha_j^S \right),$ | $j = 1, \dots, q$ |
| | $\mathbf{Q}_{\boldsymbol{\theta}} = \operatorname{cov}[\mathbf{v}_{\boldsymbol{\theta}} \mathbf{v}'_{\boldsymbol{\theta}}]$ where $(\mathbf{Q}_{\boldsymbol{\theta}})_{i,j} = \alpha_{i,j}^S$ | |

1.3. *Proposed statistical eigen-inference techniques.* In this article our focus is on developing population eigenvalue estimation and testing algorithms for models of the form in (1.2) that are robust to high-dimensionality, sample size constraints and population eigenvector misspecification in the spirit of the initial exploratory work in Rao and Edelman (2006). We are able to develop such computationally feasible algorithms by exploiting the properties of the eigenvalues of large Wishart matrices. These results analytically describe the nonrandom *blurring* of the sample eigenvalues, relative to the population eigenvalues, in the $p, n(p) \to \infty$ limit while compensating for the random *fluctuations* about the limiting behavior due to finite-dimensionality effects. The initial work in Rao and Edelman (2006) only exploited the nonrandom blurring of the sample eigenvalues without accounting for the random fluctuations; this was equivalent to employing the estimation procedure in Table 1 with $\mathbf{Q}_{\boldsymbol{\theta}} = \mathbf{I}$.

Taking into account the statistics of the fluctuations results in an improved performance and allows us to handle the situation where the sample eigenvalues are blurred to the point that the block subspace structure of the population eigenvalues cannot be visually discerned, as in Figure 1(a), thereby extending the "signal" detection capability beyond the special cases tackled in Silverstein and Combettes (1992). The nature of the mathematics being exploited makes them robust to the high-dimensionality and sample size constraints while the reliance on the sample eigenvalues alone makes them insensitive to any assumptions on the population eigenvectors. In such situations where the eigenvectors are accurately modeled, the practitioner may use the proposed methodologies to complement and "robustify" the inference provided by estimation and testing methodologies that exploit the eigenvector structure.

We consider testing the hypothesis for the equality of the population eigenvalues and statistical inference about the population eigenvalues. In other words, for some unknown $\mathbf{U}$, if $\boldsymbol{\Sigma}_0 = \mathbf{U} \boldsymbol{\Lambda}_{\boldsymbol{\theta}_0} \mathbf{U}'$, where $\boldsymbol{\Lambda}_{\boldsymbol{\theta}}$ is modeled as in (1.2), techniques to (1) test if $\boldsymbol{\Sigma} = \boldsymbol{\Sigma}_0$, and (2) estimate $\boldsymbol{\theta}_0$ are summarized in Table 1. We note that inference on the population eigenvalues is performed using the *entire sample eigen-spectrum* unlike (1.6) and (1.9). This reflects the inherent nonlinearities of the sample eigenvalue blurring induced by high-dimensionality and sample size constraints.

TABLE 2
*Comparison of performance of different techniques for estimating the nonunity population eigenvalue in Figure 1 when the block structure is known a priori*

| $p$ | $n$ | **Known U Max Like.** | **Known U Max Like.$\times p^2$** | **Anderson** | **Unknown U Max Like.** | **SEI** | **SEI$\times p^2$** |
|---|---|---|---|---|---|---|---|
| | | | (a) Bias | | | | |
| 10 | 20 | 0.0117 | 0.1168 | $-1.9994$ | $-0.5811$ | $-0.0331$ | $-0.3308$ |
| 20 | 40 | 0 | 0.0001 | $-1.9994$ | $-0.5159$ | $-0.0112$ | $-0.2244$ |
| 40 | 80 | 0.0008 | 0.0301 | $-1.9994$ | $-0.5245$ | $-0.0019$ | $-0.0776$ |
| 80 | 160 | $-0.0003$ | $-0.0259$ | $-1.9994$ | $-0.4894$ | $-0.0003$ | $-0.0221$ |
| 160 | 320 | 0.0000 | 0.0035 | $-1.9994$ | $-0.4916$ | $-0.0003$ | $-0.0411$ |
| 320 | 640 | 0.0001 | 0.0426 | $-1.9994$ | $-0.5015$ | 0.0001 | 0.0179 |
| | | | (b) MSE | | | | |
| 10 | 20 | 0.0380 | 3.7976 | 3.9990 | 0.3595 | 0.0495 | 4.9463 |
| 20 | 40 | 0.0100 | 3.9908 | 3.9990 | 0.2722 | 0.0126 | 5.0256 |
| 40 | 80 | 0.0025 | 3.9256 | 3.9991 | 0.2765 | 0.0030 | 4.8483 |
| 80 | 160 | 0.0006 | 4.1118 | 3.9991 | 0.2399 | 0.0008 | 5.1794 |
| 160 | 320 | 0.0002 | 4.1022 | 3.9990 | 0.2417 | 0.0002 | 5.0480 |
| 320 | 640 | 0.0000 | 4.0104 | 3.9990 | 0.2515 | 0.0000 | 5.0210 |

Table 2 compares the bias and mean square error of various techniques of estimating the nonunity population eigenvalue in Figure 1 (the SCM is complex-valued) when the block structure is known a priori, that is, when $t_1 = t_2 = 0.5$, and $a_2 = 1$ are known and $a := a_1$ is unknown and to be estimated. The first two columns refer to the procedure in (1.8) where the correct population eigenvectors $\mathbf{U} \neq \mathbf{I}$ are used, the third column refers to Anderson's procedure in (1.6) while the fourth column refers to the procedure in (1.8) where $\mathbf{U} = \mathbf{I}$ is used instead of the population eigenvectors. The last two columns refer to the proposed statistical eigen-inference (SEI) technique in Table 1 with $\theta := a$, $v(\theta) = \mathrm{Tr}\, S - p(0.5a + 0.5)$, and $Q_\theta = (1/2a^2 + 1/2a^2c + ac + 1/2 + 1/2c - a)c^2$ where $c = p/n$. Note that though the SEI techniques do not exploit any eigenvector information, their performance compares favorably to the maximum-likelihood technique that does. As for the other techniques, it is evident that the inherent finite sample biases in the problem cripple the estimators derived on the basis of classical large sample asymptotics.

An important implication of this in practice is that in high-dimensional, sample size starved settings, local inference, performed on a subset of sample eigenvalues alone, that fails to take into account the global structure (i.e., by modeling the remaining eigenvalues) is likely to be inaccurate, or worse misleading. In such settings, practitioners are advised to consider tests (such as the ones proposed) for the equality of the entire population eigen-spectrum instead of testing for the equality of individual population eigenvalues.

We view the inference techniques developed herein as the first step in the development of improved high-dimensional covariance matrix estimation algorithms. The issue of inverse covariance matrix estimation which Srivastava (2007) examines in the context of discriminant analysis is also related.

The approach we have in mind differs from the (sample eigenvalue) shrinkage-based techniques in Haff (1980), Dey and Srinivasan (1985) in a crucial regard. Our perspective is that the eigenvalues and the eigenvectors (or subspaces) of the sample covariance matrices are blurred relative to the population eigenvalues and eigenvectors (or subspaces), respectively. For the model considered in this article, the precise analytical characterization of the blurring of the eigenvalues (Theorem 2.7) allows us to formulate and solve the deblurring problem. The tools from free probability are applied in the first author's dissertation [see Nadakuditi (2007)] to precisely describe the blurring of the population eigenvectors (or subspaces) as well. The answer is encoded in the form of a conditional eigenvector "distribution" that explicitly takes into account the dimensionality of the system and the sample size available—the conditioning is with respect to the population eigenvalues. The idea that the covariance matrix estimate thus constructed from the deblurred eigenvalues and eigenvectors should be significantly better has merit. The development of computationally realizable eigenvector deblurring algorithms is a significant obstacle to progress along this direction of research.

1.4. *Related work.* There are other alternatives found in the literature to the block subspace hypothesis testing problem considered in this article. El Karoui (2007) provides a test for the largest eigenvalue for a large class of complex Wishart matrices including those with a population covariance matrix of the form in (1.2). Though the results are stated for the case when $p < n$, simulations confirm the validity of the techniques to the alternative case when $p \geq n$ and for real Wishart matrices. El Karoui's tests can be classified as a local test that utilizes global information, that is, information about the entire (assumed) population eigen-spectrum. Testing is performed by computing the largest eigenvalue of the sample covariance matrix, recentering, rescaling it and rejecting the hypothesis if it is too large. The recentering and rescaling parameters are determined by the $a_i$ and $t_i$ values in (1.2) while the threshold is determined by the quantiles of the appropriate (real or complex) Tracy–Widom distribution [Tracy and Widom (1994, 1996), Johnstone (2001)]. A disadvantage of this procedure is the great likelihood that recentering by the false parameter pushes the test statistic toward the left tail of the distribution. Consequently, the identity covariance hypothesis will be accepted with great likelihood whenever the recentering and rescaling coefficients are calculated for the model in (1.2) with $a_i > 1$. The proposed global test based on global information avoids this pitfall and is based on distributional results for the traces of powers of Wishart matrices that also appear in Srivastava (2005). The issue of whether a local test or a global test is more powerful is important and highlighted using simulations in the context of a joint estimation and testing problem in Section 7; its full resolution is beyond the scope of this article.

Silverstein and Combettes (1992) consider the situation when the sample eigenvalues discernibly split into distinct clusters and suggest that the proportion of the eigenvalues in each cluster will provide a good estimate of the parameters $a_i$ in (1.2). The nature of the distributional results in Bai and Silverstein (1998) imply that whenever the sample eigenvalues are thus clustered, then for large enough $p$, the estimate of $a_i$ thus obtained will be exactly equal to true value. Such a procedure could not, however, be applied for situations such as those depicted in Figure 1(a) where the sample eigenvalues do not separate into clusters. Silverstein and Combettes (1992) do not provide a strategy for computing the $t_i$ in (1.2) once the $a_i$ is computed—the proposed techniques fill the void.

A semiparametric, grid-based technique for inferring the empirical distribution function of the population eigenvalues from the sample eigenspectrum was proposed by El Karoui (2006). The procedure described can be invaluable to the practitioner in the initial data exploration stage by providing a good estimate of the number of blocks in (1.2) and a less refined estimate of the underlying $a_i$ and $t_i$ associated with each block. Our techniques can then be used to improve or test the estimates.

1.5. *Outline.* The rest of this article is organized as follows. In Section 2 we introduce the necessary definitions and summarize the relevant theorem. Concrete algorithms for computing the analytic expectations that appear in the algorithms summarized in Table 1 are presented in Section 3. The eigen-inference techniques are developed in Section 4. The performance of the algorithms is illustrated using Monte Carlo simulations in Section 5. Some concluding remarks are presented in Section 8.

## 2. Preliminaries.

DEFINITION 2.1. Let $\mathbf{A} \equiv \mathbf{A}_N$ be an $N \times N$ matrix with real eigenvalues. The $j$th sample moment is defined as

$$\operatorname{tr}(\mathbf{A}^j) := \frac{1}{N} \operatorname{Tr}(\mathbf{A}^j),$$

where Tr is the usual unnormalized trace.

DEFINITION 2.2. Let $\mathbf{A} \equiv \mathbf{A}_N$ be a sequence of self-adjoint $N \times N$ random matrices. If the limit of all moments defined as

$$\alpha_j^A =: \lim_{N \to \infty} \mathrm{E}[\mathrm{tr}(\mathbf{A}_N^j)] \qquad (N \in \mathbb{N})$$

exists, then we say that $\mathbf{A}$ *has a limit eigenvalue distribution*.

NOTATION 2.3. For a random matrix $\mathbf{A}$ with a limit eigenvalue distribution we denote by $M_A(x)$ the moment power series, which we define by

$$M_A(x) := 1 + \sum_{j \geq 1} \alpha_j^A x^j.$$

NOTATION 2.4. For a random matrix ensemble $\mathbf{A}$ with limit eigenvalue distribution we denote by $g_A(x)$ the corresponding Cauchy transform, which we define as formal power series by

$$g_A(x) := \lim_{N \to \infty} \mathrm{E}\left[\frac{1}{N}\mathrm{Tr}(x\mathbf{I}_N - \mathbf{A}_N)^{-1}\right] = \frac{1}{x}M_A(1/x).$$

DEFINITION 2.5. Let $\mathbf{A} := \mathbf{A}_N$ be an $N \times N$ self-adjoint random matrix ensemble. We say that it has a *second-order limit distribution* if for all $i, j \in \mathbb{N}$ the limits

$$\alpha_j^A := \lim_{N \to \infty} k_1(\mathrm{tr}(\mathbf{A}_N^j))$$

and

$$\alpha_{i,j}^A := \lim_{N \to \infty} k_2(\mathrm{Tr}(\mathbf{A}_N^i), \mathrm{Tr}(\mathbf{A}_N^j))$$

exist and if

$$\lim_{N \to \infty} k_r(\mathrm{Tr}(\mathbf{A}_N^{j(1)}), \ldots, \mathrm{Tr}(\mathbf{A}_N^{j(r)})) = 0$$

for all $r \geq 3$ and all $j(1), \ldots, j(r) \in \mathbb{N}$. In this definition, we denote the (classical) cumulants by $k_n$. Note that $k_1$ is just the expectation, and $k_2$ the covariance.

NOTATION 2.6. When $\mathbf{A} \equiv \mathbf{A}_N$ has a limit eigenvalue distribution, then the limits $\alpha_j^A := \lim_{N \to \infty} \mathrm{E}[\mathrm{tr}(\mathbf{A}_N^j)]$ exist. When $\mathbf{A}_N$ has a second-order limit distribution, the fluctuation

$$\mathrm{tr}(\mathbf{A}_N^j) - \alpha_j^A$$

is asymptotically Gaussian of order $1/N$. We consider the second-order covariances defined as

$$\alpha_{i,j}^A := \lim_{N \to \infty} \mathrm{cov}(\mathrm{Tr}(\mathbf{A}_N^i), \mathrm{Tr}(\mathbf{A}_N^j)),$$

and denote by $\mathcal{M}_A(x, y)$ the second-order moment power series, which we define by

$$\mathcal{M}_A(x, y) := \sum_{i,j \geq 1} \alpha_{i,j}^A x^i y^j.$$

THEOREM 2.7. *Assume that the $p \times p$ (nonrandom) covariance matrix $\boldsymbol{\Sigma} = (\boldsymbol{\Sigma}_p)_{p \in \mathbb{N}}$ has a limit eigenvalue distribution. Let $\mathbf{S}$ be the (real or complex) sample covariance matrix formed from the $n$ samples as in (1.1). Then for $p, n \to \infty$ with $p/n \to c \in (0, \infty)$, $\mathbf{S}$ has both a limit eigenvalue distribution and a second-order limit distribution. The Cauchy transform of the limit eigenvalue distribution, $g(x) \equiv g_S(x)$, satisfies the equation*

$$(2.1) \qquad g(x) = \frac{1}{1 - c + cxg(x)} g_\Sigma\left(\frac{x}{1 - c + cxg(x)}\right),$$

*with the corresponding power series $M_S(x) = 1/x g_S(1/x)$. Define $\widetilde{\mathbf{S}} = \frac{1}{n}\mathbf{X}'\mathbf{X}$ so that its moment power series is given by*

$$(2.2) \qquad M_{\widetilde{S}}(y) = c(M_S(z) - 1) + 1.$$

*The second-order moment generating series is given by*

$$(2.3a) \qquad \mathcal{M}_S(x, y) = \mathcal{M}_{\widetilde{S}}(x, y) = \frac{2}{\beta}\mathcal{M}_S^\infty(x, y),$$

*where*

$$(2.3b) \qquad \mathcal{M}_S^\infty(x, y) = xy\left(\frac{\frac{d}{dx}(xM_{\widetilde{S}}(x)) \cdot \frac{d}{dy}(yM_{\widetilde{S}}(y))}{(xM_{\widetilde{S}}(x) - yM_{\widetilde{S}}(y))^2} - \frac{1}{(x - y)^2}\right),$$

*where $\beta$ equals 1 (or 2) when the elements of $S$ are real (or complex).*

PROOF. Theorem 2.7 is due to Bai and Silverstein. They stated and proved it in Bai and Silverstein (2004) by complex analysis tools. (Note, however, that there is a missing factor 2 in their formula (2.3a) that has been corrected in their book [Bai and Silverstein (2006), page 251, Lemma 9.11, (9.8.4)].) □

Our equivalent formulation in terms of formal power series can, for the case $\beta = 2$, also be derived quite canonically by using the theory of second-order freeness. Let us also mention that the proof using second-order freeness extends easily to the situation where $\boldsymbol{\Sigma}$ is itself a random matrix with a second-order limit distribution. If we denote by $\mathcal{M}_\Sigma$ the corresponding second-order moment power series of $\Sigma$, as in Notation 2.6, then the theory of second-order freeness gives (for $\beta = 2$) the following extension of formula (2.3b):

$$
\begin{aligned}
\mathcal{M}_S^\infty(x, y) = {} & M_\Sigma(x M_{\tilde{\mathbf{S}}}(x), y M_{\tilde{\mathbf{S}}}(y)) \cdot \frac{\frac{d}{dx}(x M_{\tilde{\mathbf{S}}}(x))}{M_{\tilde{\mathbf{S}}}(x)} \cdot \frac{\frac{d}{dy}(y M_{\tilde{\mathbf{S}}}(y))}{M_{\tilde{\mathbf{S}}}(y)} \\
& + xy \left( \frac{\frac{d}{dx}(x M_{\tilde{\mathbf{S}}}(x)) \cdot \frac{d}{dy}(y M_{\tilde{\mathbf{S}}}(y))}{(x M_{\tilde{\mathbf{S}}}(x) - y M_{\tilde{\mathbf{S}}}(y))^2} - \frac{1}{(x - y)^2} \right).
\end{aligned}
$$

(2.4)

Whereas from an analytic point of view, formulas (2.1) and (2.4) might look quite mysterious, from the perspective of free probability theory there is an easy conceptual way of looking on them. Namely, they are just rewritings into formal power series of the following fact: the matrix $\tilde{\mathbf{S}}$ has a compound free Poisson distribution, for both its moments and its fluctuations. This means that the free cumulants of $\tilde{\mathbf{S}}$ of first and second-order are, up to scaling, given by the moments and the fluctuations, respectively, of $\boldsymbol{\Sigma}$. (This should be compared to: a classical compound Poisson distribution is characterized by the fact that its classical cumulants are a multiple of the moments of the corresponding "jump distribution.") In the case where $\boldsymbol{\Sigma}$ is nonrandom the fluctuations of $\boldsymbol{\Sigma}$ are clearly zero (and thus the second-order free cumulants of $\tilde{\mathbf{S}}$ vanish), that is, $\mathcal{M}_\Sigma = 0$, resulting in the special case (2.3b) of formula (2.4).

For the definition of free cumulants and more information on second-order freeness, the interested reader should consult Collins et al. (2007), in particular, Section 2.

## 3. Computational aspects.

PROPOSITION 3.1. *For* $\boldsymbol{\Sigma}_\theta = \mathbf{U} \boldsymbol{\Lambda}_\theta \mathbf{U}'$ *as in* (1.2), *let* $\boldsymbol{\theta} = (t_1, \ldots, t_{k-1}, a_1, \ldots, a_k)$, *where* $t_i = p_i / p$. *Then* $\mathbf{S}$ *has a limit eigenvalue distribution as well as a second-order limit distribution. The moments* $\alpha_j^S$, *and hence* $\alpha_{i,j}^S$, *depend on* $\boldsymbol{\theta}$ *and* $c$. *Let* $\mathbf{v}_\theta$ *be a* $q$-*by*-$1$ *vector whose* $j$th *element is given by*

$$
(\mathbf{v}_\theta)_j = \operatorname{Tr} \mathbf{S}^j - p \alpha_j^S.
$$

*Then for large* $p$ *and* $n$,

(3.1) $$ \mathbf{v}_\theta \sim \mathcal{N}(\boldsymbol{\mu}_\theta, \mathbf{Q}_\theta), $$

*where* $\boldsymbol{\mu}_\theta = 0$ *if* $\mathbf{S}$ *is complex and* $(\mathbf{Q}_\theta)_{i,j} = \alpha_{i,j}^S$.

PROOF. This follows directly from Theorem 2.7. From (3.2) and (3.4), the moments $\alpha_k^S$ depend on $\alpha^\Sigma$ and $c = p/n$ and hence on the unknown parameter vector $\boldsymbol{\theta}$. The existence of the nonzero mean when $\mathbf{S}$ is real follows from the statement in Bai and Silverstein (2004). □

3.1. *Computation of moments of limiting eigenvalue distribution.* A method of enumerating the moments of the limiting eigenvalue distribution is to use the software package RMTool [Rao (2006)] based on the polynomial method developed in the second part of the first author's dissertation [Nadakuditi (2007)]. The software enables the moments of $\mathbf{S}$ to be enumerated rapidly whenever the moment power series of $\boldsymbol{\Sigma}$ is an algebraic power series, that is, it is the solution of an algebraic equation. This is always the case when $\boldsymbol{\Sigma}$ is of the form in (1.2). For example, if $\boldsymbol{\theta} = (t_1, t_2, a_1, a_2, a_3)$, then we can obtain the moments of $\mathbf{S}$ by typing in the following sequence of commands in Matlab once RMTool has been installed. This eliminates the need to obtain manually the expressions for the moments a priori:

```
>> startRMTool
>> syms c t1 t2 a1 a2 a3
>> number_of_moments = 5;
>> LmzSigma = atomLmz([a1 a2 a3],[t1 t2 1-(t1+t2)]);
>> LmzS = AtimesWish(LmzSigma,c);
>> alpha_S = Lmz2MomF(LmzS,number_of_moments);
>> alpha_Stilde = c*alpha_S;
```

An alternate and versatile method of computing the moments relies on exploiting (2.1) which expresses the relationship between the moment power series of $\boldsymbol{\Sigma}$ and that of $\mathbf{S}$ via the limit of the ratio $p/n$. This allows us to directly express the expected moments of $\mathbf{S}$ in terms of the moments of $\boldsymbol{\Sigma}$. The general form of the moments of $\widetilde{\mathbf{S}}$, given by Corollary 9.12 in Nica and Speicher [(2006), page 143] is

$$(3.2) \quad \alpha_j^{\widetilde{S}} = \sum_{\substack{i_j \geq 0 \\ 1i_1 + 2i_2 + 3i_3 + \cdots + ji_j = j}} c^{i_1 + i_2 + \cdots + i_j} (\alpha_1^{\Sigma})^{i_1} (\alpha_2^{\Sigma})^{i_2} \cdots (\alpha_j^{\Sigma})^{i_j} \cdot \gamma_{i_1, i_2, \ldots, i_j}^{(j)},$$

where $\gamma_{i_1, \ldots, i_j}^{j}$ is the multinomial coefficient given by

$$(3.3) \qquad \gamma_{i_1, i_2, \ldots, i_j}^{(j)} = \frac{j!}{i_1! i_2! \cdots i_j! (j + 1 - (i_1 + i_2 + \cdots + i_j))!}.$$

The multinomial coefficient in (3.3) has an interesting combinatorial interpretation. Let $j$ be a positive integer, and let $i_1, \ldots, i_j \in \mathbb{N} \cup \{0\}$ be such that $i_1 + 2i_2 + \cdots + ji_j = j$. The number of noncrossing partitions $\pi \in NC(j)$ which have $i_1$ blocks with 1 element, $i_2$ blocks with 2 elements, ..., $i_j$ blocks with $j$ elements is given by the multinomial coefficient $\gamma_{i_1, \ldots, i_j}^{j}$.

The moments of $\widetilde{\mathbf{S}}$ are related to the moments of $\mathbf{S}$ as

$$(3.4) \qquad\qquad \alpha_j^{\widetilde{S}} = c\alpha_j^S \qquad \text{for } j = 1, 2, \ldots.$$

We can use (3.2) to compute the first few moments of $\mathbf{S}$ in terms of the moments of $\boldsymbol{\Sigma}$. This involves enumerating the partitions that appear in the computation of the multinomial coefficient in (3.3). For $j = 1$ only $i_1 = 1$ contributes with $\gamma_1^{(1)} = 1$, thus,

$$(3.5) \qquad \alpha_1^{\widetilde{S}} = c\alpha_1^{\Sigma}.$$

For $n = 2$ only $i_1 = 2$, $i_2 = 0$ and $i_1 = 0$, $i_2 = 1$ contribute with

$$\gamma_{2,0}^{(2)} = 1, \qquad \gamma_{0,1}^{(2)} = 1$$

and thus

$$(3.6) \qquad \alpha_2^{\widetilde{S}} = c\alpha_2^{\Sigma} + c^2(\alpha_1^{\Sigma})^2.$$

For $n = 3$ we have three possibilities for the indices, contributing with

$$\gamma_{3,0,0}^{(3)} = 1, \qquad \gamma_{1,1,0}^{(3)} = 3, \qquad \gamma_{0,0,1}^{(3)} = 1,$$

thus

$$(3.7) \qquad \alpha_3^{\widetilde{S}} = c\alpha_3^{\Sigma} + 3c^2\alpha_1^{\Sigma}\alpha_2^{\Sigma} + c^3(\alpha_1^{\Sigma})^3.$$

For $n = 4$ we have five possibilities for the indices, contributing with

$$\gamma_{4,0,0,0}^{(4)} = 1, \qquad \gamma_{2,1,0,0}^{(4)} = 6, \qquad \gamma_{0,2,0,0}^{(4)} = 2, \qquad \gamma_{1,0,1,0}^{(4)} = 4, \qquad \gamma_{0,0,0,1}^{(4)} = 1,$$

thus

$$(3.8) \qquad \alpha_4^{\widetilde{S}} = c\alpha_4^{\Sigma} + 4c^2\alpha_1^{\Sigma}\alpha_3^{\Sigma} + 2c^2(\alpha_2^{\Sigma})^2 + 6c^3(\alpha_1^{\Sigma})^2\alpha_2^{\Sigma} + c^4(\alpha_1^{\Sigma})^4.$$

For specific instances of $\boldsymbol{\Sigma}$, we simply plug the moments $\alpha_i^{\Sigma}$ into the above expressions to get the corresponding moments of $\mathbf{S}$. The general formula in (3.2) can be used to generate the expressions for higher-order moments. We can efficiently enumerate the sum-constrained partitions that appear in (3.2) by employing the algorithm that recursively computes the nonnegative integer sequences $s(k)$ of length $j$ with the sum constraint $\sum_{k=1}^{j} s(k) k = j$ listed below

    Input: Integer $j$

    Output: Nonnegative integer sequences $s(k)$ of length $j$ satisfying constraint $\sum_k [s(k) \times k] = n$

If $j = 1$

    The only sequence of length 1 is $s = j$

else

    for $k = 0$ to 1

        Compute sequences of length $j - 1$ for $j - k \times j$

        Append $k$ to each sequence above and include in output

    end

end

3.2. *Computation of covariance moments of second-order limit distribution.* Equations (2.3a) and (2.3b) express the relationship between the covariance of the second-order limit distribution and the moments of **S**. Let $M(x)$ denote a moment power series as in Notation 2.3 with coefficients $\alpha_j$. Define the power series $H(x) = xM(x)$ and let

$$(3.9) \qquad \mathcal{H}(x,y) := \left( \frac{\frac{d}{dx}(H(x)) \cdot \frac{d}{dy}(H(y))}{(H(x) - H(y))^2} - \frac{1}{(x-y)^2} \right)$$

so that $\mathcal{M}^\infty(x,y) := xy\mathcal{H}(x,y)$. The $(i,j)$th coefficient of $\mathcal{M}^\infty(x,y)$ can then be extracted from a multivariate Taylor series expansion of $\mathcal{H}(x,y)$ about $x=0$, $y=0$. From (2.3a), we then obtain the coefficients $\alpha_{i,j}^S = (2/\beta)\alpha_{i,j}^{\mathcal{M}^\infty}$. The coefficients $\alpha_{i,j}^S$ can be readily enumerated by invoking a short sequence of commands in the MAPLE computer algebra system. For example, the code on the next page will enumerate $\alpha_{5,2}^S$. By modifying this code, we can obtain the coefficients $\alpha_{i,j}^S$ in terms of $\alpha_i := \alpha_i^S = \alpha_j^{\widetilde{S}}/c$ for other choices of indices $i$ and $j$ and the constant max_coeff chosen such that $i+j \leq 2\max\_\text{coeff}$.

```
> with(numapprox):
> max_coeff := 5:
> H := x -> x*(1+sum(alpha[j]*x^2,j=1..2*max_coeff)):
> dHx : = diff(H(x),x): dHy := diff(H(y),y):
> H2 := simplify(dHx*dHy/(H(x)-H(y))^2-1/(x-y)^2:
> H2series := mtaylor(H2,[x,y],2*max_coeff):
> i:=5: j :=2:
> M2_infty_coeff[i,j]
   := simplify(coeff(coeff(H2series,x,i-1),y,j-1)):
> alphaS_second[i,j] := (2/beta)*M2_infty_coeff[i,j]:
```

Table 3 lists some of the coefficients of $\mathcal{M}^\infty$ obtained using this procedure. When $\alpha_j = 1$ for all $j \in \mathbb{N}$, then $\alpha_{i,j} = 0$ as expected, since $\alpha_j = 1$ denotes the identity matrix. Note that the moments $\alpha_1, \ldots, \alpha_{i+j}$ are needed to compute the second-order covariance moments $\alpha_{i,j} = \alpha_{j,i}$.

The covariance matrix **Q** with elements $\mathbf{Q}_{i,j} = \alpha_{i,j}$ gets increasingly ill-conditioned as $\dim(\mathbf{Q})$ increases; the growth in the magnitude of the diagonal entries $\alpha_{j,j}$ in Table 3 attests to this. This implies that the eigenvectors of **Q** encode the information about the covariance of the second-order limit distribution more efficiently than the matrix **Q** itself. When $\boldsymbol{\Sigma} = \mathbf{I}$ so that the SCM **S** has the (null) Wishart distribution, the eigenvectors of **Q** are the (appropriately normalized) Chebyshev polynomials of the second kind [Mingo and Speicher (2006)]. The structure of the eigenvectors for arbitrary $\boldsymbol{\Sigma}$ is, as yet, unknown though research in that direction might yield additional insights.

TABLE 3
*Relationship between the coefficients* $\alpha_{i,j} = \alpha_{j,i}$ *and* $\alpha_i$

| Coefficient | Expression |
|---|---|
| $\boldsymbol{\alpha_{1,1}}$ | $\alpha_2 - \alpha_1^2$ |
| $\boldsymbol{\alpha_{2,1}}$ | $-4\alpha_1\alpha_2 + 2\alpha_1^3 + 2\alpha_3$ |
| $\boldsymbol{\alpha_{2,2}}$ | $16\alpha_1^2\alpha_2 - 6\alpha_2^2 - 6\alpha_1^4 - 8\alpha_1\alpha_3 + 4\alpha_4$ |
| $\boldsymbol{\alpha_{3,1}}$ | $9\alpha_1^2\alpha_2 - 6\alpha_1\alpha_3 - 3\alpha_2^2 + 3\alpha_4 - 3\alpha_1^4$ |
| $\boldsymbol{\alpha_{3,2}}$ | $6\alpha_5 + 30\alpha_1\alpha_2^2 - 42\alpha_1^3\alpha_2 - 18\alpha_2\alpha_3 + 12\alpha_1^5 + 24\alpha_1^2\alpha_3 - 12\alpha_1\alpha_4$ |
| $\boldsymbol{\alpha_{3,3}}$ | $-18\alpha_2^3 - 27\alpha_2\alpha_4 + 9\alpha_6 - 30\alpha_1^6 + 21\alpha_3^2 + 36\alpha_1^2\alpha_4 - 72\alpha_1^3\alpha_3 + 126\alpha_1^4\alpha_2 - 135\alpha_1^2\alpha_2^2 + 108\alpha_1\alpha_2\alpha_3 - 18\alpha_1\alpha_5$ |
| $\boldsymbol{\alpha_{4,1}}$ | $12\alpha_1\alpha_2^2 - 16\alpha_1^3\alpha_2 - 8\alpha_2\alpha_3 + 12\alpha_1^2\alpha_3 - 8\alpha_1\alpha_4 + 4\alpha_1^5 + 4\alpha_5$ |
| $\boldsymbol{\alpha_{4,2}}$ | $-12\alpha_2^3 - 24\alpha_2\alpha_4 + 8\alpha_6 - 20\alpha_1^6 + 16\alpha_3^2 + 32\alpha_1^2\alpha_4 - 56\alpha_1^3\alpha_3 + 88\alpha_1^4\alpha_2 - 96\alpha_1^2\alpha_2^2 + 80\alpha_1\alpha_2\alpha_3 - 16\alpha_1\alpha_5$ |
| $\boldsymbol{\alpha_{4,3}}$ | $96\alpha_2^2\alpha_3 + 60\alpha_1^7 + 84\alpha_1\alpha_3^2 + 432\alpha_1^3\alpha_2^2 + 180\alpha_1^4\alpha_3 - 48\alpha_3\alpha_4 + 12\alpha_7 - 36\alpha_2\alpha_5 - 24\alpha_1\alpha_6 + 144\alpha_1\alpha_2\alpha_4 + 48\alpha_1^2\alpha_5 - 96\alpha_1^3\alpha_4 - 156\alpha_1\alpha_2^3 - 300\alpha_1^5\alpha_2 - 396\alpha_1^2\alpha_2\alpha_3$ |
| $\boldsymbol{\alpha_{4,4}}$ | $-140\alpha_1^8 - 76\alpha_2^4 - 48\alpha_6\alpha_2 + 256\alpha_3\alpha_4\alpha_1 - 40\alpha_4^2 + 16\alpha_8 - 64\alpha_3\alpha_5 - 32\alpha_1\alpha_7 + 1408\alpha_1^3\alpha_2\alpha_3 - 336\alpha_1^2\alpha_3^2 + 256\alpha_1^4\alpha_4 + 144\alpha_2^2\alpha_4 - 480\alpha_1^5\alpha_3 + 160\alpha_2\alpha_3^2 + 64\alpha_1^2\alpha_6 - 128\alpha_1^3\alpha_5 - 1440\alpha_1^4\alpha_2^2 + 832\alpha_1^2\alpha_2^3 + 800\alpha_1^6\alpha_2 - 768\alpha_1\alpha_2^2\alpha_3 - 576\alpha_1^2\alpha_2\alpha_4 + 192\alpha_1\alpha_2\alpha_5$ |
| $\boldsymbol{\alpha_{5,1}}$ | $-5\alpha_3^2 - 10\alpha_2\alpha_4 + 5\alpha_6 - 5\alpha_1^6 + 5\alpha_2^3 + 15\alpha_1^2\alpha_4 - 20\alpha_1^3\alpha_3 + 25\alpha_1^4\alpha_2 - 30\alpha_1^2\alpha_2^2 + 30\alpha_1\alpha_2\alpha_3 - 10\alpha_1\alpha_5$ |
| $\boldsymbol{\alpha_{5,2}}$ | $60\alpha_2^2\alpha_3 + 30\alpha_1^7 + 50\alpha_1\alpha_3^2 + 240\alpha_1^3\alpha_2^2 + 110\alpha_1^4\alpha_3 - 30\alpha_3\alpha_4 + 10\alpha_7 - 30\alpha_2\alpha_5 - 20\alpha_1\alpha_6 + 100\alpha_1\alpha_2\alpha_4 + 40\alpha_1^2\alpha_5 - 70\alpha_1^3\alpha_4 - 90\alpha_1\alpha_2^3 - 160\alpha_1^5\alpha_2 - 240\alpha_1^2\alpha_2\alpha_3$ |
| $\boldsymbol{\alpha_{5,3}}$ | $-105\alpha_1^8 - 60\alpha_2^4 - 45\alpha_6\alpha_2 + 210\alpha_3\alpha_4\alpha_1 - 30\alpha_4^2 + 15\alpha_8 - 60\alpha_3\alpha_5 - 30\alpha_1\alpha_7 + 1140\alpha_1^3\alpha_2\alpha_3 - 270\alpha_1^2\alpha_3^2 + 225\alpha_1^4\alpha_4 + 120\alpha_2^2\alpha_4 - 390\alpha_1^5\alpha_3 + 135\alpha_2\alpha_3^2 + 60\alpha_1^2\alpha_6 - 120\alpha_1^3\alpha_5 - 1125\alpha_1^4\alpha_2^2 + 660\alpha_1^2\alpha_2^3 + 615\alpha_1^6\alpha_2 - 630\alpha_1\alpha_2^2\alpha_3 - 495\alpha_1^2\alpha_2\alpha_4 + 180\alpha_1\alpha_2\alpha_5$ |
| $\boldsymbol{\alpha_{5,4}}$ | $-900\alpha_1^2\alpha_4\alpha_3 + 80\alpha_1^2\alpha_7 - 160\alpha_1^3\alpha_6 - 620\alpha_1^5\alpha_4 - 3200\alpha_1^3\alpha_2^3 + 700\alpha_1\alpha_2^4 + 3960\alpha_1^5\alpha_2^2 - 720\alpha_1^2\alpha_5\alpha_2 + 1840\alpha_1^3\alpha_4\alpha_2 - 4100\alpha_1^4\alpha_3\alpha_2 + 3600\alpha_1^2\alpha_2^2\alpha_3 - 1140\alpha_1\alpha_2^3\alpha_2 + 1040\alpha_1^3\alpha_3^2 - 440\alpha_2^3\alpha_3 + 440\alpha_3\alpha_4\alpha_2 + 240\alpha_1\alpha_6\alpha_2 + 320\alpha_1\alpha_5\alpha_3 - 1020\alpha_1\alpha_2^2\alpha_4 + 20\alpha_9 - 1820\alpha_1^7\alpha_2 + 180\alpha_2^2\alpha_5 + 320\alpha_1^4\alpha_5 + 180\alpha_1\alpha_4^2 + 1120\alpha_1^6\alpha_3 + 80\alpha_3^3 + 280\alpha_1^9 - 40\alpha_1\alpha_8 - 60\alpha_7\alpha_2 - 80\alpha_3\alpha_6 - 100\alpha_4\alpha_5$ |
| $\boldsymbol{\alpha_{5,5}}$ | $2400\alpha_2\alpha_5\alpha_3^2 - 1350\alpha_2^2\alpha_5\alpha_1 + 600\alpha_3\alpha_5\alpha_2 + 300\alpha_1\alpha_7\alpha_2 - 900\alpha_6\alpha_2\alpha_1^2 - 1200\alpha_3\alpha_5\alpha_1^2 + 400\alpha_1\alpha_6\alpha_3 + 3000\alpha_3\alpha_4\alpha_1^3 + 5100\alpha_1^2\alpha_2^2\alpha_4 + 12300\alpha_1^5\alpha_2\alpha_3 + 5700\alpha_1^2\alpha_2^2\alpha_3^2 + 4400\alpha_1\alpha_2^3\alpha_3 + 400\alpha_1^4\alpha_6 - 15000\alpha_1^3\alpha_2^2\alpha_3 - 5750\alpha_1^4\alpha_2\alpha_4 - 200\alpha_1^3\alpha_7 + 500\alpha_1\alpha_4\alpha_5 + 225\alpha_6\alpha_2^2 - 675\alpha_2^2\alpha_1^2 - 3250\alpha_1^4\alpha_3^2 - 625\alpha_2^2\alpha_4 + 350\alpha_3^2\alpha_4 - 600\alpha_1\alpha_3^3 - 1050\alpha_2^2\alpha_3^2 - 2800\alpha_3\alpha_1^7 - 11550\alpha_1^6\alpha_2^2 - 3300\alpha_3\alpha_4\alpha_1\alpha_2 - 800\alpha_3\alpha_5^3 + 325\alpha_3^2\alpha_2 - 4375\alpha_1^2\alpha_2^4 - 630\alpha_1^{10} + 100\alpha_8\alpha_1^2 - 75\alpha_5^2 + 255\alpha_2^5 + 12000\alpha_1^4\alpha_2^3 + 4550\alpha_1^8\alpha_2 + 1550\alpha_1^6\alpha_4 + 25\alpha_{10} - 50\alpha_1\alpha_9 - 75\alpha_2\alpha_8 - 100\alpha_3\alpha_7 - 125\alpha_4\alpha_6$ |

## 4. Eigen-inference algorithms.

4.1. *Estimating $\theta$ for known model order.* Estimating the unknown parameter vector $\boldsymbol{\theta}$ follows from the asymptotic result in Proposition 3.1. For large $p, n$, since $\mathbf{v}_{\boldsymbol{\theta}}$ is (approximately) normally distributed we can obtain the estimate $\boldsymbol{\theta}$ by the principle of maximum-likelihood. When $\mathbf{S}$ is real, Bai and Silverstein provide a formula, expressed as a difficult-to-compute contour integral, for the correction term $\boldsymbol{\mu}_{\boldsymbol{\theta}}$ in (3.1). The log-likelihood of $\mathbf{v}_{\boldsymbol{\theta}}$ is (ignoring constants and the correction term for the mean when $\mathbf{S}$ is real) given by

$$(4.1) \qquad \ell(\mathbf{v}_{\boldsymbol{\theta}}|\boldsymbol{\theta}) \approx -\mathbf{v}_{\boldsymbol{\theta}}^T \mathbf{Q}_{\boldsymbol{\theta}}^{-1} \mathbf{v}_{\boldsymbol{\theta}} - \log \det \mathbf{Q}_{\boldsymbol{\theta}},$$

which allows us to obtain the maximum-likelihood estimate of $\boldsymbol{\theta}$ as

$$(4.2) \quad \widehat{\boldsymbol{\theta}}_{(q)} = \underset{\boldsymbol{\theta} \in \boldsymbol{\Theta}}{\arg\min}\, \mathbf{v}_{\boldsymbol{\theta}}^T \mathbf{Q}_{\boldsymbol{\theta}}^{-1} \mathbf{v}_{\boldsymbol{\theta}} + \log \det \mathbf{Q}_{\boldsymbol{\theta}} \qquad \text{for } q = \dim(\mathbf{v}_{\boldsymbol{\theta}}) \geq \dim(\boldsymbol{\theta}),$$

where $\boldsymbol{\Theta}$ represents the parameter space for the elements of $\boldsymbol{\theta}$ and $\mathbf{v}_{\boldsymbol{\theta}}$ and $\mathbf{Q}_{\boldsymbol{\theta}}$ are constructed as in Proposition 3.1.

4.2. *Guidelines for picking $q := \dim(\mathbf{v}_{\theta})$.* Canonically, the parameter vector $\boldsymbol{\theta}$ of models such as (1.2) is of length $2k - 1$ so that $q = \dim(\mathbf{v}_{\boldsymbol{\theta}})$ must equal or exceed $2k - 1$. In principle, estimation accuracy should increase with $q$ as the covariance of $\mathbf{v}_{\boldsymbol{\theta}}$ is explicitly accounted for via the weighting matrix $\mathbf{Q}_{\boldsymbol{\theta}}$.

Figure 2 compares the quantiles of the test statistic $\mathbf{v}_{\boldsymbol{\theta}}' \mathbf{Q}_{\boldsymbol{\theta}}^{-1} \mathbf{v}_{\boldsymbol{\theta}}$ for $\dim(\mathbf{v}_{\boldsymbol{\theta}}) = q$ with the quantiles of the chi-squared distribution with $q$ degrees of freedom when $q = 2, 3$ for the model in (1.2) with $\boldsymbol{\theta} = (0.5, 2, 1)$, $n = p$ and $p = 40$ and $p = 320$. While there is good agreement with the theoretical distribution for large $n, p$, the deviation from the limiting result is not insignificant for moderate $n, p$. This justifies setting $q = 2$ for the testing procedures developed herein.

In the most general estimation setting as in (4.2) where $\boldsymbol{\theta}$ includes the smallest population eigenvalue of $\boldsymbol{\Sigma}_{\boldsymbol{\theta}}$ we have found that $q := \dim(\mathbf{v}_{\boldsymbol{\theta}})$ must be no smaller than $\dim(\boldsymbol{\theta}) + 1$. When the smallest eigenvalue of $\boldsymbol{\Sigma}_{\boldsymbol{\theta}}$ is known, however, $q$ can be as small as $\dim(\boldsymbol{\theta})$. Within these guidelines, picking a smaller value of $q$ provides robustness in low-to moderate-dimensional settings where the deviations from the asymptotic result in Theorem 2.7 are not insignificant. Numerical simulations suggest that the resulting degradation in estimation accuracy in high-dimensional settings, in using the smallest suggested choice for $q$ instead of a higher value, is relatively small. This loss in performance is offset by an increase in the speed of the underlying numerical optimization routine. This is the case because, though the dimensionality of $\boldsymbol{\theta}$ is the same, the matrix $\mathbf{Q}$ gets increasingly ill-conditioned for higher values of $q$, thereby reducing the efficiency of optimization methods.

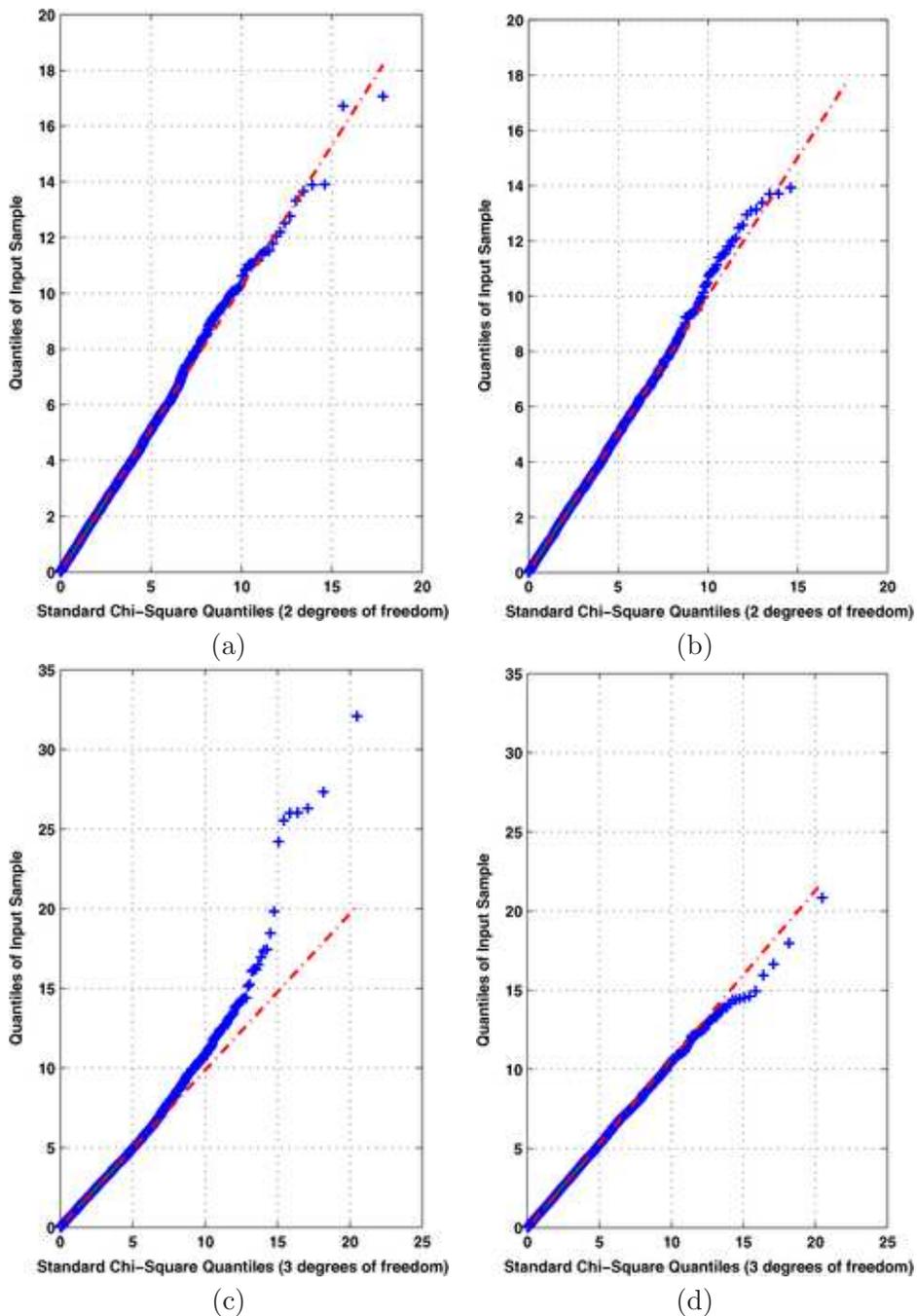

Fig. 2. *Numerical simulations (when* **S** *is complex) illustrating the robustness of the distribution approximation for the test statistic in Table 1 formed with* $\dim(\mathbf{v}) = 2$ *to moderate-dimensional settings.* (a) $p = n = 40$: $\dim(\mathbf{v}) = 2$. (b) $p = n = 320$: $\dim(\mathbf{v}) = 2$. (c) $p = n = 40$: $\dim(\mathbf{v}) = 3$. (d) $p = n = 320$: $\dim(\mathbf{v}) = 3$.

4.3. *Testing $\theta = \theta_0$.*

PROPOSITION 4.1. *Define the vector $\mathbf{v_\theta}$ and the covariance matrix $\mathbf{Q_\theta}$ as*

$$(4.3a) \quad \mathbf{v_\theta} = \begin{bmatrix} \mathrm{Tr}\,\mathbf{S} - p\alpha_1^\Sigma \\ \mathrm{Tr}\,\mathbf{S}^2 - p\left(\alpha_2^\Sigma + \dfrac{p}{n}(\alpha_1^\Sigma)^2\right) - \left(\dfrac{2}{\beta} - 1\right)\alpha_2^\Sigma\dfrac{p}{n} \end{bmatrix},$$

$$(4.3b) \quad \mathbf{Q_\theta} = \frac{2}{\beta}\begin{bmatrix} \widetilde{\alpha}_2 - \widetilde{\alpha}_1^2 & 2\widetilde{\alpha}_1^3 + 2\widetilde{\alpha}_3 - 4\widetilde{\alpha}_1\widetilde{\alpha}_2 \\ 2\widetilde{\alpha}_1^3 + 2\widetilde{\alpha}_3 - 4\widetilde{\alpha}_1\widetilde{\alpha}_2 & 4\widetilde{\alpha}_4 - 8\widetilde{\alpha}_1\widetilde{\alpha}_3 - 6\widetilde{\alpha}_2^2 + 16\widetilde{\alpha}_2\widetilde{\alpha}_1^2 - 6\widetilde{\alpha}_1^4 \end{bmatrix},$$

*with $\beta = 1$ (or 2) when S is real (or complex) and $\widetilde{\alpha}_i \equiv \alpha_i^{\widetilde{S}}$ given by*

$$(4.4a) \quad \widetilde{\alpha}_1 = \frac{p}{n}\alpha_1^\Sigma,$$

$$(4.4b) \quad \widetilde{\alpha}_2 = \frac{p}{n}\alpha_2^\Sigma + \frac{p^2}{n^2}(\alpha_1^\Sigma)^2,$$

$$(4.4c) \quad \widetilde{\alpha}_3 = \frac{p}{n}\alpha_3^\Sigma + 3\frac{p^2}{n^2}\alpha_1^\Sigma\alpha_2^\Sigma + \frac{p^3}{n^3}(\alpha_1^\Sigma)^3,$$

$$(4.4d) \quad \widetilde{\alpha}_4 = \frac{p}{n}\alpha_4^\Sigma + 4\frac{p^2}{n^2}\alpha_1^\Sigma\alpha_3^\Sigma + 2\frac{p^2}{n^2}(\alpha_2^\Sigma)^2 + 6\frac{p^3}{n^3}(\alpha_1^\Sigma)^2\alpha_2^\Sigma + \frac{p^4}{n^4}(\alpha_1^\Sigma)^4,$$

*and $\alpha_i^\Sigma = (1/p)\,\mathrm{Tr}\,\mathbf{\Sigma}^i = \sum_{j=1}^{k} a_j t_j^i$. Thus, for large $p$ and $n$, $\mathbf{v_\theta} \sim \mathcal{N}(\mathbf{0}, \mathbf{Q_\theta})$ so that*

$$(4.5) \qquad h(\boldsymbol{\theta}) := \mathbf{v_\theta}^T \mathbf{Q_\theta}^{-1} \mathbf{v_\theta} \sim \chi_2^2.$$

PROOF. This follows from Proposition 3.1. The correction term for the real case is discussed in a different context in Dumitriu, Edelman and Shuman (2007). A matrix-theoretic derivation in the real case ($\beta = 1$) can be found in Srivastava (2005), Corollary 2.1, page 3. □

We test for $\boldsymbol{\theta} = \boldsymbol{\theta}_0$ by obtaining the test statistic

$$(4.6) \qquad H_{\boldsymbol{\theta}_0} : h(\boldsymbol{\theta}_0) = \mathbf{v}_{\boldsymbol{\theta}_0}^T \mathbf{Q}_{\boldsymbol{\theta}_0}^{-1} \mathbf{v}_{\boldsymbol{\theta}_0},$$

where the $\mathbf{v}_{\boldsymbol{\theta}_0}$ and $\mathbf{Q}_{\boldsymbol{\theta}_0}$ are constructed as in (4.3a) and (4.3b), respectively. We reject the hypothesis for large values of $H_{\boldsymbol{\theta}_0}$. For a choice of threshold $\gamma$, the asymptotic convergence of the test statistic to the $\chi_2^2$ distribution implies that

$$(4.7) \qquad \mathrm{Prob.}(H_{\boldsymbol{\theta}_0} = 1 | \boldsymbol{\theta} = \boldsymbol{\theta}_0) \approx F^{\chi_2^2}(\gamma).$$

Thus, for large $p$ and $n$, when $\gamma = 5.9914$, $\mathrm{Prob.}(H_{\boldsymbol{\theta}_0} = 1 | \boldsymbol{\theta} = \boldsymbol{\theta}_0) \approx 0.95$.

4.4. *Estimating $\theta$ and testing the estimate.* When a $\widehat{\boldsymbol{\theta}}$ is obtained using (4.2) then we may test for $\theta = \widehat{\boldsymbol{\theta}}$ by forming the testing statistic

$$(4.8) \qquad H_{\widehat{\boldsymbol{\theta}}}\colon h(\widehat{\boldsymbol{\theta}}) = \mathbf{u}_{\widehat{\boldsymbol{\theta}}}^T \mathbf{W}_{\widehat{\boldsymbol{\theta}}}^{-1} \mathbf{u}_{\widehat{\boldsymbol{\theta}}},$$

where the $\mathbf{u}_{\widehat{\boldsymbol{\theta}}}$ and $\mathbf{W}_{\widehat{\boldsymbol{\theta}}}$ are constructed as in (4.3a) and (4.3b), respectively. However, the sample covariance matrix $\mathbf{S}$ can no longer be used since the estimate $\widehat{\boldsymbol{\theta}}$ was obtained from it. Instead, we form a test sample covariance matrix constructed from $\lceil (n/2) \rceil$ randomly chosen samples. Equivalently, since the samples are assumed to be mutually independent and identically distributed, we can form the test matrix from the first $\lceil (n/2) \rceil$ samples as

$$(4.9) \qquad \overline{\mathbf{S}} = \frac{1}{\lceil n/2 \rceil} \sum_{i=1}^{\lceil n/2 \rceil} \mathbf{x}_i \mathbf{x}_i'.$$

Note that $\alpha_k^{\overline{\mathbf{S}}}$ will have to be recomputed using $\boldsymbol{\Sigma}_{\widehat{\boldsymbol{\theta}}}$ and $\overline{c} = p / \lceil (n/2) \rceil$. The hypothesis $\boldsymbol{\theta} = \widehat{\boldsymbol{\theta}}$ is tested by rejecting values of the test statistic greater than a threshold $\gamma$. The threshold is selected using the approximation in (4.7). Alternately, the hypothesis can be rejected if the recentered and rescaled largest eigenvalue of $\mathbf{S}$ is greater than the threshold $\gamma$. The threshold is selected using the quantiles of the (real or complex) Tracy–Widom distribution. The recentering and rescaling coefficients are obtained by the procedure described in El Karoui (2007).

4.5. *Estimating $\theta$ for unknown model order.* Suppose we have a family of models parameterized by the vector $\boldsymbol{\theta}^{(\overline{k})}$. The elements of $\boldsymbol{\theta}^{(\overline{k})}$ are the free parameters of the model. For the model in (1.2), in the canonical case $\boldsymbol{\theta} = (t_1, \ldots, t_{k-1}, a_1, \ldots, a_k)$ since $t_1 + \cdots + t_{k-1} + t_k = 1$ so that $\dim(\boldsymbol{\theta}^{(\overline{k})}) = 2k - 1$. If some of the parameters in (1.2) are known, then the parameter vector is modified accordingly.

When the model order is unknown, we select the model which has the minimum Akaike Information Criterion. For the situation at hand we propose that

$$(4.10) \qquad \begin{aligned} & \widehat{\boldsymbol{\theta}} = \widehat{\boldsymbol{\theta}}^{(\widehat{k})} \\ & \text{where } \widehat{k} = \operatorname*{arg\,min}_{k \in \mathbb{N}} \{ \mathbf{u}_{\widehat{\boldsymbol{\theta}}^{(k)}}^T \mathbf{W}_{\widehat{\boldsymbol{\theta}}^{(k)}}^{-1} \mathbf{u}_{\widehat{\boldsymbol{\theta}}^{(k)}} + \log \det \mathbf{W}_{\widehat{\boldsymbol{\theta}}^{(k)}} \} + 2 \dim(\boldsymbol{\theta}^{(k)}), \end{aligned}$$

where $\mathbf{u}_{\widehat{\boldsymbol{\theta}}^{(k)}}$ and $\mathbf{W}_{\widehat{\boldsymbol{\theta}}^{(k)}}$ are constructed as described in Section 4.4 using the test sample covariance matrix in (4.9). Alternately, a sequence of nested hypothesis tests using a largest eigenvalue based test as described in El Karoui (2007) can be used. It would be useful to compare the performance of the proposed and the nested hypothesis testing procedures in situations of practical interest.

In the point of view adopted in this article, the sample eigen-spectrum is a *single observation* sampled from the multivariate probability distribution in (1.5). Thus we did not consider a Bayesian Information Criterion (BIC) based formulation in (4.10) because of the resulting degeneracy of the conventional "log(Sample Size)" penalty term. In the context of model selection, the study of penalty function selection, including issues that arise due to dimensionality, remains an important topic whose full resolution is beyond the scope of this article. Nevertheless, we are able to demonstrate the robustness of the method proposed in (4.10) in some representative situations in Section 6.

**5. Numerical simulations.** Let $\Sigma_{\overline{\theta}}$ be as in (1.2) with $\overline{\theta} = (t_1, a_1, a_2)$. When $t_1 = 0.5$, $a_1 = 2$ and $a_2 = 1$, then half of the population eigenvalues are equal to 2 while the remainder are of magnitude 1. Let the unknown parameter vector $\theta = (t, a)$ where $t \equiv t_1$ and $a \equiv a_1$. Using the procedure described in Section 3.1, the first four moments can be obtained as (here $c = p/n$)

$$\text{(5.1a)} \qquad \alpha_1^S = 1 + t(a - 1),$$

$$\text{(5.1b)} \qquad \alpha_2^S = (-2ac + a^2c + c)t^2 + (-1 + 2ac - 2c + a^2)t + 1 + c,$$

$$\alpha_3^S = (-3c^2a^2 + a^3c^2 - c^2 + 3ac^2)t^3$$

$$\text{(5.1c)} \qquad + (3c^2 + 3c^2a^2 - 3ac - 6ac^2 - 3a^2c + 3a^3c + 3c)t^2$$

$$+ (-3c^2 + a^3 - 1 - 6c + 3ac + 3a^2c + 3ac^2)t + 1 + c^2 + 3c,$$

$$\alpha_4^S = (6a^2c^3 + a^4c^3 - 4ac^3 - 4a^3c^3 + c^3)t^4$$

$$+ (-6c^2 - 12a^3c^2 + 12ac^3$$

$$\qquad - 12a^2c^3 + 4a^3c^3 + 12ac^2 + 6a^4c^2 - 4c^3)t^3$$

$$+ (-4a^2c - 4ac - 12ac^3 - 24ac^2 + 6a^4c$$

$$\text{(5.1d)} \qquad + 6a^2c^3 + 12a^3c^2 + 6c - 6c^2a^2 + 6c^3 + 18c^2 - 4a^3c)t^2$$

$$+ (-4c^3 + 4ac + 6c^2a^2 + 4ac^3$$

$$\qquad - 1 + 12ac^2 - 18c^2 + 4a^2c - 12c + 4a^3c + a^4)t$$

$$+ 1 + c^3 + 6c + 6c^2.$$

From the discussion in Section 3.2, we obtain the covariance of the second-order limit distribution

$$\text{(5.2)} \quad \mathbf{Q}_{\theta} = \frac{2}{\beta} \begin{bmatrix} c^2(\alpha_2^S - \alpha_1^2) & c^3(2(\alpha_1^S)^3 + 2\alpha_3^S - 4\alpha_1^S\alpha_2^S) \\ c^3(2(\alpha_1^S)^3 + 2\alpha_3^S - 4\alpha_1^S\alpha_2^S) & c^4(4\alpha_4^S - 8\alpha_1^S\alpha_3^S \\ & \quad -6(\alpha_2^S)^2 + 16\alpha_2^S(\alpha_1^S)^2 - 6(\alpha_1^S)^4) \end{bmatrix},$$

TABLE 4
*Quality of estimation of $t = 0.5$ for different values of $p$ (dimension of observation vector) and $n$ (number of samples)—both real and complex case for the example in Section 5*

| | | Complex case | | | Real case | | |
|---|---|---|---|---|---|---|---|
| $p$ | $n$ | Bias | MSE | $(MSE \times p^2)/100$ | Bias | MSE | $(MSE \times p^2)/100$ |
| | | | | (a) $n = 0.5p$ | | | |
| 20 | 10 | 0.0455 | 0.3658 | 1.4632 | 0.4862 | 1.2479 | 4.9915 |
| 40 | 20 | −0.0046 | 0.1167 | 1.8671 | 0.2430 | 0.3205 | 5.1272 |
| 80 | 40 | −0.0122 | 0.0337 | 2.1595 | 0.1137 | 0.08495 | 5.437 |
| 160 | 80 | −0.0024 | 0.0083 | 2.1250 | 0.0598 | 0.02084 | 5.335 |
| 320 | 160 | 0.0008 | 0.0021 | 2.1790 | 0.0300 | 0.00528 | 5.406 |
| | | | | (b) $n = p$ | | | |
| 20 | 20 | −0.0137 | 0.1299 | 0.5196 | 0.2243 | 0.3483 | 1.3932 |
| 40 | 40 | −0.0052 | 0.0390 | 0.6233 | 0.1083 | 0.0901 | 1.4412 |
| 80 | 80 | −0.0019 | 0.0093 | 0.5941 | 0.0605 | 0.0231 | 1.4787 |
| 160 | 160 | −0.0005 | 0.0024 | 0.6127 | 0.0303 | 0.0055 | 1.4106 |
| 320 | 320 | −0.0001 | 0.0006 | 0.6113 | 0.0162 | 0.0015 | 1.5155 |
| | | | | (c) $n = 2p$ | | | |
| 20 | 40 | −0.0119 | 0.0420 | 0.1679 | 0.1085 | 0.1020 | 0.4081 |
| 40 | 80 | −0.0017 | 0.0109 | 0.1740 | 0.0563 | 0.0255 | 0.4079 |
| 80 | 160 | −0.0005 | 0.0028 | 0.1765 | 0.0290 | 0.0063 | 0.4056 |
| 160 | 320 | −0.0004 | 0.0007 | 0.1828 | 0.0151 | 0.0016 | 0.4139 |
| 320 | 640 | 0.0001 | 0.0002 | 0.1752 | 0.0080 | 0.0004 | 0.4024 |

where $\beta = 1$ when $S$ is real-valued and $\beta = 2$ when **S** is complex-valued.

We then use (4.2) to estimate $\boldsymbol{\theta}$ and hence the unknown parameters $t$ and $a$. Tables 4 and 5 compare the bias and mean squared error of the estimates for $a$ and $t$, respectively. Note the $1/p^2$ type decay in the mean squared error and how the real case has twice the variance as the complex case. As expected by the theory of maximum-likelihood estimation, the estimates become increasingly normal for large $p$ and $n$. This is evident from Figure 3. As expected, the performance improves as the dimensionality of the system increases.

## 6. Model order related issues.

6.1. *Robustness to model order overspecification.* Consider the situation when the samples are complex-valued and the true covariance matrix $\boldsymbol{\Sigma} = 2\mathbf{I}$. We erroneously assume that there are two blocks for the model in (1.2) and that $a_2 = 1$ is known while $a := a_1$ and $t := t_1$ are unknown and have to be estimated. We estimate $\boldsymbol{\theta} = (a, t)$ using (4.2) as before. The empirical cumulative distribution function (CDF) of $\hat{t}$ over 4000 Monte

TABLE 5
*Quality of estimation of $a = 2$ for different values of $p$ (dimension of observation vector) and $n$ (number of samples)—both real and complex case for the example in Section 5*

| | | Complex case | | | Real case | | |
|---|---|---|---|---|---|---|---|
| $p$ | $n$ | Bias | MSE | $(\text{MSE} \times p^2)/100$ | Bias | MSE | $(\text{MSE} \times p^2)/100$ |
| | | | | (a) $n = 0.5p$ | | | |
| 20 | 10 | 0.1278 | 0.1046 | 0.4185 | 0.00748 | 0.1024 | 0.4097 |
| 40 | 20 | 0.0674 | 0.0478 | 0.7647 | −0.01835 | 0.04993 | 0.7989 |
| 80 | 40 | 0.0238 | 0.0111 | 0.7116 | −0.02240 | 0.01800 | 1.1545 |
| 160 | 80 | 0.0055 | 0.0022 | 0.5639 | −0.02146 | 0.00414 | 1.0563 |
| 320 | 160 | 0.0007 | 0.0005 | 0.5418 | −0.01263 | 0.00112 | 1.1692 |
| | | | | (b) $n = p$ | | | |
| 20 | 20 | 0.0750 | 0.0525 | 0.2099 | −0.0019 | 0.0577 | 0.2307 |
| 40 | 40 | 0.0227 | 0.0127 | 0.2028 | −0.0206 | 0.0187 | 0.2992 |
| 80 | 80 | 0.0052 | 0.0024 | 0.1544 | −0.0206 | 0.0047 | 0.3007 |
| 160 | 160 | 0.0014 | 0.0006 | 0.1499 | −0.0126 | 0.0012 | 0.3065 |
| 320 | 320 | 0.0003 | 0.0001 | 0.1447 | −0.0074 | 0.0003 | 0.3407 |
| | | | | (c) $n = 2p$ | | | |
| 20 | 40 | 0.0251 | 0.0134 | 0.0534 | −0.0182 | 0.0205 | 0.0821 |
| 40 | 80 | 0.0049 | 0.0028 | 0.0447 | −0.0175 | 0.0052 | 0.0834 |
| 80 | 160 | 0.0015 | 0.0007 | 0.0428 | −0.0115 | 0.0014 | 0.0865 |
| 160 | 320 | 0.0004 | 0.0002 | 0.0434 | −0.0067 | 0.0004 | 0.0920 |
| 320 | 640 | 0.0000 | 0.0000 | 0.0412 | −0.0038 | 0.0001 | 0.0932 |

Carlo trials shown in Figure 4(d) shows that $\hat{t} \to 1$ as $p, n(p) \to \infty$. Figure 4(c) compares the quantiles of test statistic in (4.5) with that of the chi-squared distribution with two degrees of freedom. The excellent agreement for modest values of $p$ and $n$ validates the distributional approximation. Figure 4(a) and (b) plot the mean squared errors in estimating $a$ and $t$, respectively. As before, the mean squared error exhibits a $1/p^2$ behavior. Table 6 shows the $1/p$ decay in the bias of estimating these parameters.

For this same example, the seventh and eighth columns of Table 6 show the level at which a sphericity and the 2 block hypothesis are accepted when the procedure described in (4.2) is applied and a threshold is set at the 5% significance level. The ninth and tenth columns of Table 6 show the acceptance rate for the 2 block hypothesis when the largest eigenvalue test proposed in El Karoui (2007) is applied on a test sample covariance matrix formed using first $\lceil n/2 \rceil$ samples and the original sample covariance matrix, respectively. The largest eigenvalue value test has an acceptance rate closer to the 5% significance level it was designed for. For all of the $p$ and $n$ values in Table 6, over the 4000 Monte Carlo trials, applying the procedure described in Section 4.5 produced the correct estimate $\hat{k} = 1$ for the order of the model in (1.2) when $\boldsymbol{\Sigma} = 2\mathbf{I}$.

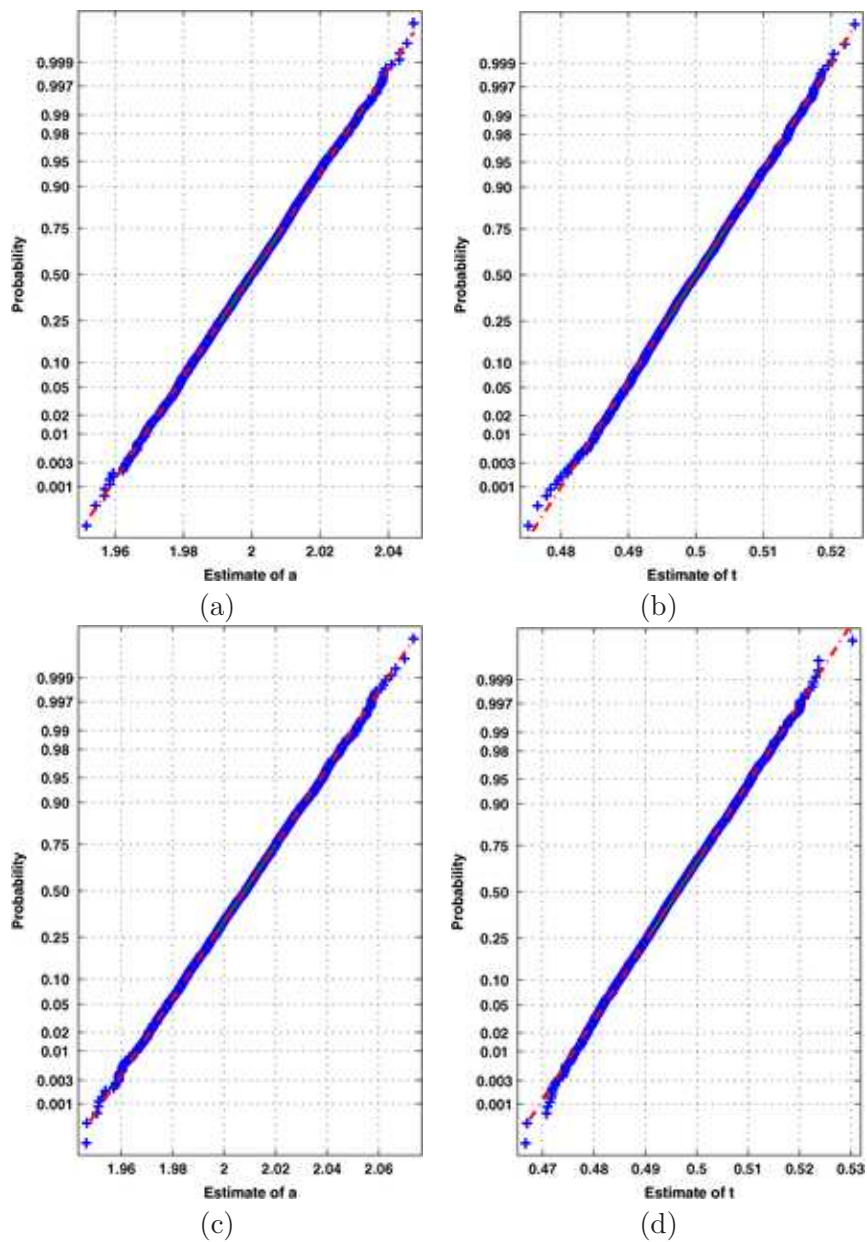

Fig. 3. *Normal probability plots of the estimates of a and t (true values: $a = 2$, $t = 0.5$) for the example in Section 5.* (a) $\widehat{a}$: $p = 320, n = 640$. (b) $\widehat{t}$: $p = 320, p = 640$. (c) $\widehat{a}$: $p = 320, n = 640$ *(real-valued)*. (d) $\widehat{t}$: $p = 320, n = 640$ *(real-valued)*.

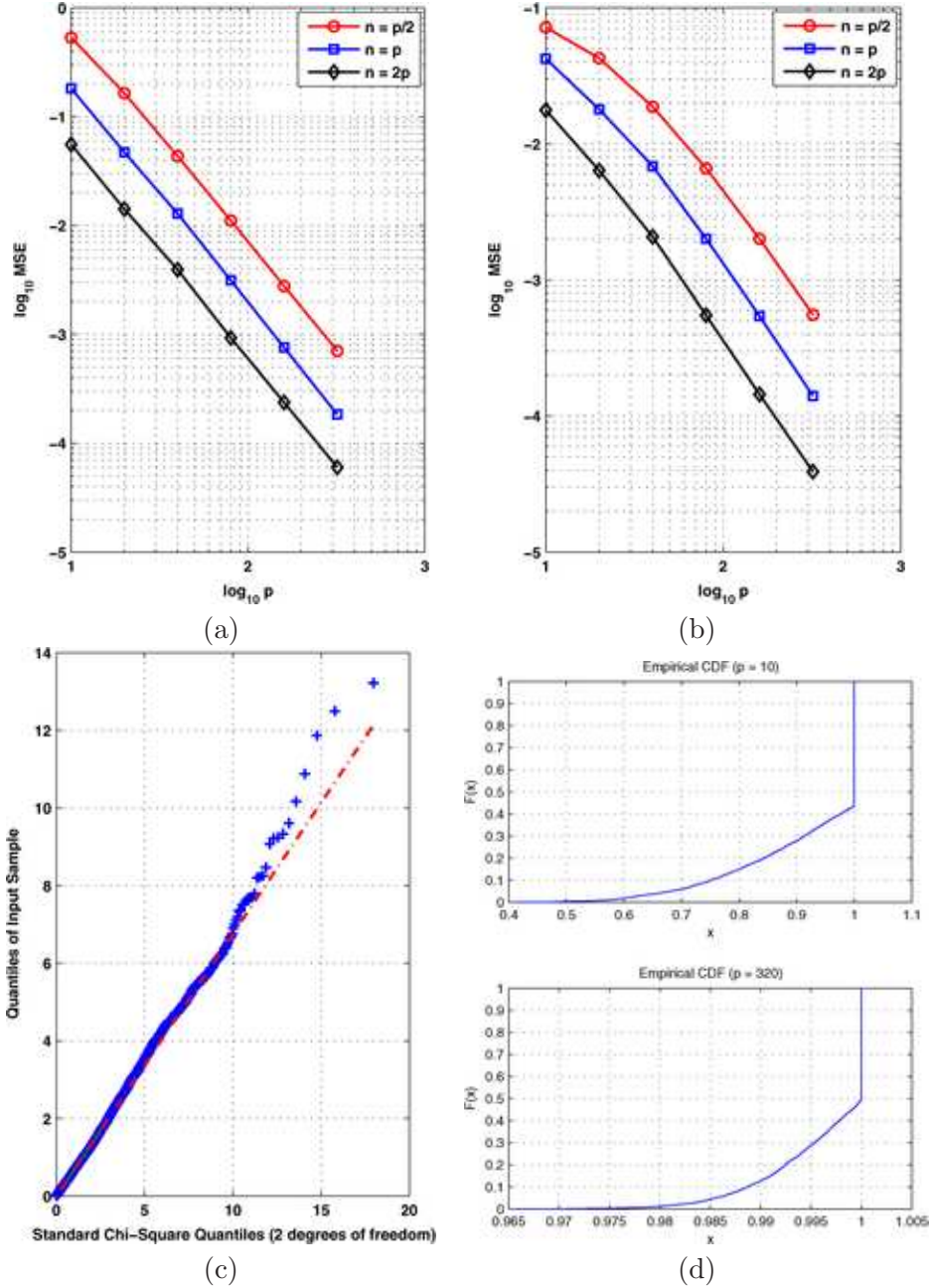

FIG. 4. *Performance of estimation algorithm when model order has been overspecified and* **S** *is complex. The population covariance matrix* $\boldsymbol{\Sigma} = 2\mathbf{I}$ *which corresponds in (1.2) to* $a_1 = 1$, *and* $t_1 = 1$ *for arbitrary* $a_2$. *We run the estimation algorithm assuming that* $a_2 = 1$ *and estimate* $a := a_1$ *and* $t := t_1$ *in (1.2).* (a) *MSE:* $\hat{a}$. (b) *MSE:* $\hat{t}$. (c) *QQ plot: Test statistic in (4.5) for* $p = 320 = 2n$. (d) *Empirical CDF of* $\hat{t}$: $n = p/2$.

Table 6

*Performance of estimation algorithm when model order has been overspecified and* **S** *is complex*

| $p$ | $n$ | $\hat{a}$ Bias | $\hat{a}$ Bias $\times$ p | $\hat{t}$ Bias | $\hat{t}$ Bias $\times$ p | Sphericity acceptance | 2 Block acceptance | $\lambda_{\max}$ test (full) | $\lambda_{\max}$ test (half) |
|---|---|---|---|---|---|---|---|---|---|
| | | | | (a) $n = p/2$ | | | | | |
| 10 | 5 | 0.3523 | 3.5232 | −0.1425 | −1.4246 | 0.9820 | 0.9801 | 1.0000 | 0.9698 |
| 20 | 10 | 0.1997 | 3.9935 | −0.1157 | −2.3148 | 0.9783 | 0.9838 | 0.9998 | 0.9710 |
| 40 | 20 | 0.1078 | 4.3114 | −0.0783 | −3.1336 | 0.9795 | 0.9870 | 0.9958 | 0.9713 |
| 80 | 40 | 0.0545 | 4.3561 | −0.0463 | −3.7018 | 0.9765 | 0.9873 | 0.9838 | 0.9720 |
| 160 | 80 | 0.0272 | 4.3530 | −0.0251 | −4.0175 | 0.9743 | 0.9828 | 0.9763 | 0.9643 |
| 320 | 160 | 0.0141 | 4.5261 | −0.0133 | −4.2580 | 0.9805 | 0.9885 | 0.9753 | 0.9675 |
| | | | | (b) $n = p$ | | | | | |
| 10 | 10 | 0.2087 | 2.0867 | −0.1123 | −1.1225 | 0.9793 | 0.9768 | 0.9998 | 0.9675 |
| 20 | 20 | 0.1050 | 2.0991 | −0.0753 | −1.5060 | 0.9773 | 0.9845 | 0.9965 | 0.9723 |
| 40 | 40 | 0.0558 | 2.2312 | −0.0470 | −1.8807 | 0.9850 | 0.9898 | 0.9898 | 0.9743 |
| 80 | 80 | 0.0283 | 2.2611 | −0.0255 | −2.0410 | 0.9813 | 0.9868 | 0.9773 | 0.9710 |
| 160 | 160 | 0.0137 | 2.1990 | −0.0130 | −2.0811 | 0.9805 | 0.9870 | 0.9790 | 0.9613 |
| 320 | 320 | 0.0067 | 2.1455 | −0.0067 | −2.1568 | 0.9775 | 0.9835 | 0.9608 | 0.9603 |
| | | | | (c) $n = 2p$ | | | | | |
| 10 | 20 | 0.1067 | 1.0674 | −0.0717 | −0.7171 | 0.9790 | 0.9810 | 0.9993 | 0.9708 |
| 20 | 40 | 0.0541 | 1.0811 | −0.0442 | −0.8830 | 0.9753 | 0.9858 | 0.9890 | 0.9708 |
| 40 | 80 | 0.0290 | 1.1581 | −0.0257 | −1.0272 | 0.9743 | 0.9845 | 0.9830 | 0.9695 |
| 80 | 160 | 0.0140 | 1.1161 | −0.0131 | −1.0497 | 0.9763 | 0.9850 | 0.9743 | 0.9658 |
| 160 | 320 | 0.0071 | 1.1302 | −0.0068 | −1.0883 | 0.9778 | 0.9830 | 0.9703 | 0.9578 |
| 320 | 640 | 0.0036 | 1.1549 | −0.0035 | −1.1237 | 0.9758 | 0.9833 | 0.9598 | 0.9608 |

The population covariance matrix $\boldsymbol{\Sigma} = 2\mathbf{I}$ which corresponds in (1.2) to $a_1 = 1$, and $t_1 = 1$ for arbitrary $a_2$. We run the estimation algorithm assuming that $a_1 = 1$ and estimate $a := a_2$ and $t := t_1$ in (1.2).

6.2. *Robust model order estimation.* We revisit the setting in Section 5, where the population parameter vector $\boldsymbol{\theta} \equiv (a_1, a_2, t_1) = (2, 1, 0.5)$ and the sample covariance matrix is formed from complex-valued data. We employ the procedure described in Section 4.5 to estimate the model order $k$ (assumed unknown) and the corresponding $2k - 1$-dimensional parameter vector $\boldsymbol{\theta}^{(k)}$. Over 1000 Monte Carlo trials, for values of $p, n$ listed in Table 7, we observe the robustness to model order overspecification as in Section 6.1. Additionally, we note that as $p, n(p) \to \infty$, the correct model order is estimated consistently. Table 7 demonstrates that, as before, the parameter vector is estimated with greater accuracy as the dimensionality of the system is increased. The parameter estimates appear to be asymptotically unbiased and normally distributed as before.

TABLE 7
*Performance of parameter estimation algorithm when model order has to be estimated as well and **S** is complex*

| $p$ | $n$ | $\mathbf{Pr}\,(\widehat{k}=1)$ | $\widehat{a}$ | $\mathbf{Pr}\,(\widehat{k}=2)$ | $\widehat{a}_1$ | $\widehat{a}_2$ | $\widehat{t}_1$ |
|---|---|---|---|---|---|---|---|
| | | | (a) $n = p/2$ | | | | |
| 20 | 10 | 0.968 | $1.4867 \pm 0.1105$ | 0.032 | $1.8784 \pm 0.7384$ | $0.8785 \pm 0.6376$ | $0.5675 \pm 0.2650$ |
| 40 | 20 | 0.940 | $1.4985 \pm 0.0567$ | 0.060 | $2.0287 \pm 0.7244$ | $0.6929 \pm 0.6010$ | $0.6041 \pm 0.3165$ |
| 80 | 40 | 0.700 | $1.4990 \pm 0.0274$ | 0.300 | $2.0692 \pm 0.4968$ | $0.7604 \pm 0.4751$ | $0.5624 \pm 0.2965$ |
| 160 | 80 | 0.199 | $1.4998 \pm 0.0142$ | 0.801 | $2.0199 \pm 0.2780$ | $0.9062 \pm 0.2841$ | $0.5311 \pm 0.2084$ |
| 320 | 160 | 0.001 | $1.4999 \pm 0.0069$ | 0.999 | $2.0089 \pm 0.1398$ | $0.9763 \pm 0.1341$ | $0.5076 \pm 0.1239$ |
| 480 | 240 | – | $1.4999 \pm 0.0046$ | 1 | $2.0004 \pm 0.0967$ | $0.9847 \pm 0.0918$ | $0.5076 \pm 0.0887$ |
| | | | (b) $n = p$ | | | | |
| 20 | 20 | 0.915 | $1.4867 \pm 0.0806$ | 0.085 | $1.9229 \pm 0.5675$ | $0.6747 \pm 0.5748$ | $0.6293 \pm 0.2962$ |
| 40 | 40 | 0.736 | $1.4987 \pm 0.0381$ | 0.264 | $1.9697 \pm 0.3719$ | $0.7685 \pm 0.4199$ | $0.5920 \pm 0.2644$ |
| 80 | 80 | 0.190 | $1.5005 \pm 0.0197$ | 0.810 | $2.0021 \pm 0.2273$ | $0.9287 \pm 0.2323$ | $0.5310 \pm 0.1856$ |
| 160 | 160 | 0.004 | $1.4997 \pm 0.0099$ | 0.996 | $1.9908 \pm 0.1108$ | $0.9771 \pm 0.0995$ | $0.5162 \pm 0.0973$ |
| 320 | 320 | – | $1.5000 \pm 0.0048$ | 1 | $2.0001 \pm 0.0548$ | $0.9960 \pm 0.0469$ | $0.5024 \pm 0.0492$ |
| 480 | 480 | – | $1.5000 \pm 0.0033$ | 1 | $2.0018 \pm 0.0363$ | $1.0002 \pm 0.0310$ | $0.4991 \pm 0.0327$ |
| | | | (c) $n = 2p$ | | | | |
| 20 | 40 | 0.743 | $1.4972 \pm 0.0556$ | 0.257 | $1.9124 \pm 0.3044$ | $0.7835 \pm 0.3756$ | $0.6087 \pm 0.2424$ |
| 40 | 80 | 0.217 | $1.5002 \pm 0.0286$ | 0.783 | $1.9707 \pm 0.1797$ | $0.9361 \pm 0.1659$ | $0.5444 \pm 0.1458$ |
| 80 | 160 | – | $1.4996 \pm 0.0139$ | 1 | $1.9925 \pm 0.0975$ | $0.9847 \pm 0.0781$ | $0.5116 \pm 0.0807$ |
| 160 | 320 | – | $1.4999 \pm 0.0071$ | 1 | $1.9975 \pm 0.0485$ | $0.9959 \pm 0.0369$ | $0.5034 \pm 0.0401$ |
| 320 | 640 | – | $1.5001 \pm 0.0035$ | 1 | $1.9994 \pm 0.0232$ | $0.9993 \pm 0.0178$ | $0.5008 \pm 0.0193$ |
| 480 | 960 | – | $1.4999 \pm 0.0024$ | 1 | $1.9998 \pm 0.0161$ | $0.9996 \pm 0.0125$ | $0.5003 \pm 0.0135$ |

The population covariance matrix has parameters $a_2 = 2$, $a_1 = 1$ and $t_1 = 0.5$ as in (1.2). The algorithm in (4.10) with $\dim(v) = 2k$ for $k = 1, 2, \ldots, 5$ was used to estimate the model order $k$ and the associated $2k - 1$-dimensional parameter vector $\boldsymbol{\theta} = (a_1, \ldots, a_k, t_1, \ldots, t_{k-1})$. Numerical results shown were computed over 1000 Monte Carlo trials.

**7. Inferential aspects of spiked covariance matrix models.** Consider covariance matrix models whose eigenvalues are of the form $\lambda_1 \geq \lambda_2 \geq \cdots \geq \lambda_k > \lambda_{k+1} = \cdots = \lambda_p = \lambda$. Such models arise when the signal occupies a $k$-dimensional subspace and the noise has covariance $\lambda \mathbf{I}$. Such models are referred to as *spiked covariance matrix models*. When $k \ll p$, then for large $p$, for $\mathbf{v}_{\boldsymbol{\theta}}$ defined as in Proposition 3.1, the matrix $\mathbf{Q}_{\boldsymbol{\theta}}$ may be constructed from the moments of the (null) Wishart distribution [Dumitriu and Rassart (2003)] instead, which are given by

$$(7.1) \qquad \alpha_k^W = \lambda^k \sum_{j=0}^{k-1} c^j \frac{1}{j+1} \binom{k}{j} \binom{k-1}{j},$$

where $c = p/n$. Thus, for $q = 2$, $\mathbf{Q}_{\boldsymbol{\theta}}$ is given by

$$(7.2) \qquad \mathbf{Q}_{\boldsymbol{\theta}} \equiv \mathbf{Q}_{\lambda} = \frac{2}{\beta} \begin{bmatrix} \lambda^2 c & 2\lambda^3(c+1)c \\ 2\lambda^3(c+1)c & 2\lambda^4(2c^2 + 5c + 2)c \end{bmatrix}.$$

This substitution is motivated by Bai and Silverstein's analysis [Bai and Silverstein (2004)] where it is shown that when $k$ is small relative to $p$, then the second-order fluctuation distribution is asymptotically independent of the "spikes." When the multiplicities of the spike are known (say 1), then we let $t_i = 1/p$ and compute the moments $\alpha_j^S$ accordingly. The estimation problem thus reduces to

$$(7.3) \qquad \widehat{\boldsymbol{\theta}} = \arg\min_{\boldsymbol{\theta} \in \Theta} \mathbf{v}_{\boldsymbol{\theta}}^T \mathbf{Q}_{\lambda}^{-1} \mathbf{v}_{\boldsymbol{\theta}} \qquad \text{with } q = \dim(\mathbf{v}_{\boldsymbol{\theta}}) = \dim(\boldsymbol{\theta}) + 1,$$

where $\lambda$ is an element of $\boldsymbol{\theta}$ when it is unknown.

Consider the problem of estimating the magnitude of the spike for the model in (1.2) with $t_1 = 1/p$, and $a_2 = 1$ known and $a_1 = 10$ unknown so that $\boldsymbol{\theta} = a \equiv a_1$. We obtain the estimate $\widehat{\boldsymbol{\theta}}$ from (7.3) with $\lambda = 1$ wherein the moments $\alpha_k^S$ given by

$$(7.4a) \qquad \alpha_1^S = \frac{-1 + a + p}{p},$$

$$(7.4b) \qquad \alpha_2^S = \frac{a^2 p - 2pc + c - 2ac + cp^2 + p^2 - p + 2pac + a^2 c}{p^2}$$

are obtained by plugging in $t = 1/p$ into (5.1).

Table 8 summarizes the estimation performance for this example. Note the $1/p$ scaling of the mean squared error and how the complex case has half the mean squared error. The estimates produced are asymptotically normal as seen in Figure 5.

Table 8

*Algorithm performance for different values of p (dimension of observation vector) and n
(number of samples)—both real and complex case*

| $p$ | $n$ | | Complex case | | | Real case | |
|---|---|---|---|---|---|---|---|
| | | Bias | MSE | MSE × p | Bias | MSE | MSE × p |
| | | | | (a) $n = p$ | | | |
| 10 | 10 | −0.5528 | 9.3312 | 93.3120 | −0.5612 | 18.4181 | 184.1808 |
| 20 | 20 | −0.2407 | 4.8444 | 96.8871 | −0.2005 | 9.6207 | 192.4143 |
| 40 | 40 | −0.1168 | 2.5352 | 101.4074 | −0.0427 | 4.9949 | 199.7965 |
| 80 | 80 | −0.0833 | 1.2419 | 99.3510 | −0.03662 | 2.4994 | 199.9565 |
| 160 | 160 | −0.0371 | 0.6318 | 101.0949 | 0.03751 | 1.2268 | 196.3018 |
| 320 | 320 | −0.0125 | 0.3186 | 101.9388 | 0.04927 | 0.6420 | 204.4711 |
| | | | | (b) $n = 1.5p$ | | | |
| 10 | 15 | −0.3343 | 6.6954 | 66.9537 | −0.3168 | 12.7099 | 127.0991 |
| 20 | 30 | −0.1781 | 3.2473 | 64.9454 | −0.1454 | 6.4439 | 128.8798 |
| 40 | 60 | −0.1126 | 1.6655 | 66.6186 | −0.08347 | 3.2470 | 129.88188 |
| 80 | 120 | −0.0565 | 0.8358 | 66.8600 | −0.02661 | 1.6381 | 131.04739 |
| 160 | 240 | −0.0287 | 0.4101 | 65.6120 | 0.02318 | 0.8534 | 136.5475 |
| 320 | 480 | −0.0135 | 0.2083 | 66.6571 | 0.02168 | 0.4352 | 139.2527 |
| | | | | (c) $n = 2p$ | | | |
| 10 | 20 | −0.2319 | 4.9049 | 49.0494 | −0.2764 | 9.6992 | 96.9922 |
| 20 | 40 | −0.1500 | 2.5033 | 50.0666 | −0.1657 | 4.6752 | 93.5043 |
| 40 | 80 | −0.0687 | 1.2094 | 48.3761 | −0.03922 | 2.5300 | 101.2007 |
| 80 | 160 | −0.0482 | 0.6214 | 49.7090 | −0.02426 | 1.2252 | 98.0234 |
| 160 | 320 | −0.0111 | 0.3160 | 50.5613 | 0.01892 | 0.6273 | 100.3799 |
| 320 | 640 | −0.0139 | 0.1580 | 50.5636 | 0.02748 | 0.3267 | 104.5465 |

7.1. *Limitations.* Consider testing for the hypothesis that $\mathbf{\Sigma} = \mathbf{I}$ using real valued observations. For the model in (1.2), which is equivalent to testing $\boldsymbol{\theta} = (1, 1)$, from the discussion in Section 4.3, we form the test statistic

(7.5) $$H_{\text{Sph.}} : h(\boldsymbol{\theta}) = \mathbf{v}_{\boldsymbol{\theta}}^T \mathbf{Q}_{\boldsymbol{\theta}}^{-1} \mathbf{v}_{\boldsymbol{\theta}},$$

where $\mathbf{Q}_{\boldsymbol{\theta}}$ is given by (7.2) with $\beta = 1$ since the observations are real valued, $\lambda = 1$ and

$$\mathbf{v}_{\boldsymbol{\theta}} = \begin{bmatrix} \text{Tr } \mathbf{S} - p \\ \text{Tr } \mathbf{S}^2 - p\left(1 + \dfrac{p}{n}\right) - \left(\dfrac{2}{\beta} - 1\right)\dfrac{p}{n} \end{bmatrix},$$

where $c = p/n$, as usual. We set a threshold $\gamma = 5.9914$ so that we accept the sphericity hypothesis whenever $h(\boldsymbol{\theta}) \leq \gamma$. This corresponds to the 95th percentile of the $\chi_2^2$ distribution. Table 9(a) demonstrates how the test is able to accept the identity covariance hypothesis when $\mathbf{\Sigma}_{\boldsymbol{\theta}} = \mathbf{I}$ at a rate close to the 5% significance level it was designed for. Table 9(b) shows the acceptance of the sphericity hypothesis when $\mathbf{\Sigma}_{\boldsymbol{\theta}} = \text{diag}(10, 1, \ldots, 1)$ instead.

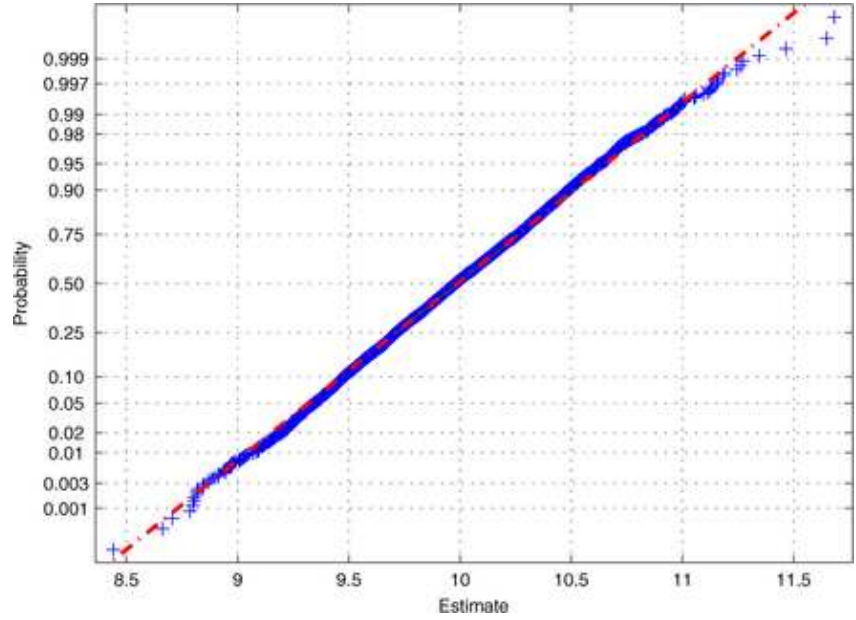

(a)

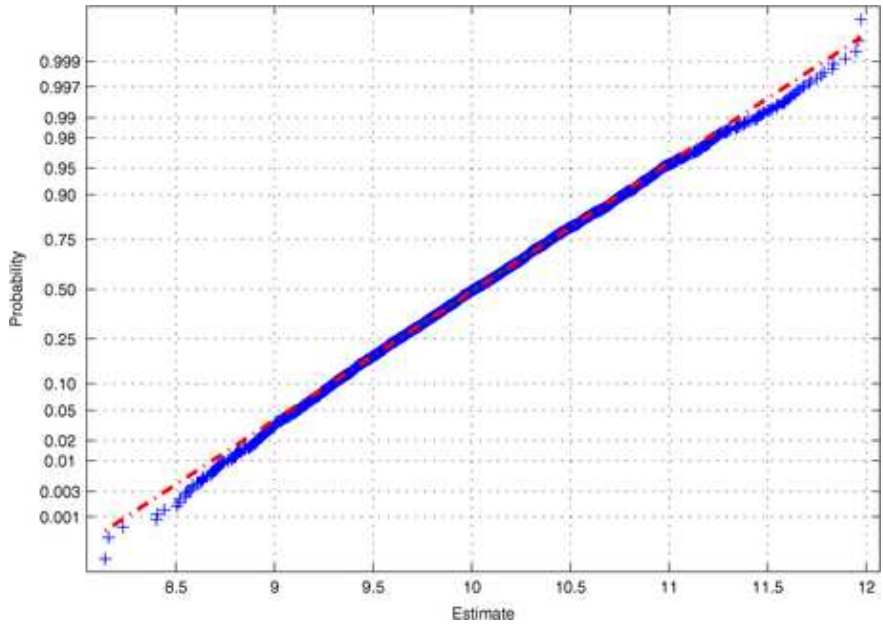

(b)

FIG. 5. *Normal probability plots of the spiked magnitude estimate (true value = 10).*
(a) $p = 320, n = 640$ *(complex* **S***).* (b) $p = 320, n = 640$ *(real* **S***).*

Results were tabulated over 4000 Monte-Carlo trials. Table 10 illustrates the performance of the sphericity test proposed by Ledoit and Wolf (2002) which consists of forming the test statistic

$$(7.6) \qquad \mathrm{LW}(\mathbf{S}) := \frac{np}{2}\left[\frac{1}{p}\operatorname{Tr}[(\mathbf{S}-\mathbf{I})^2]-\frac{p}{n}\left[\frac{1}{p}\operatorname{Tr}\mathbf{S}\right]^2+\frac{p}{n}\right] \to \chi^2_{p(p+1)/2}$$

and rejecting for large values above a threshold that is determined by using the asymptotic chi-squared approximation. Note how when $p/n$ is large, both tests erroneously accept the identity covariance hypothesis an inordinate number of times. The faulty inference provided by the test based on the methodologies developed is best understood in the context of the following result.

PROPOSITION 7.1. *Let* $\mathbf{S}$ *denote a sample covariance matrix formed from an* $p \times n$ *matrix of Gaussian observations whose columns are independent of each other and identically distributed with mean* $\mathbf{0}$ *and covariance* $\mathbf{\Sigma}$. *Denote the eigenvalues of* $\mathbf{\Sigma}$ *by* $\lambda_1 \geq \lambda_2 > \cdots \geq \lambda_k > \lambda_{k+1} = \cdots = \lambda_p = \lambda$. *Let* $l_j$ *denote the* $j$th *largest eigenvalue of* $\widehat{\mathbf{R}}$. *Then as* $p, n \to \infty$ *with* $c_n =$

TABLE 9

*The identity covariance hypothesis is rejected when the test statistic in (7.5) exceeds the 5% significance level for the $\chi^2$ distribution with 2 degress of freedom, that is, whenever $h(\theta) > 5.9914$*

|  | $n=10$ | $n=20$ | $n=40$ | $n=80$ | $n=160$ | $n=320$ | $n=640$ |
|---|---|---|---|---|---|---|---|
| (a) Empirical probability of accepting the identity covariance hypothesis when $\mathbf{\Sigma}_\theta = \mathbf{I}$ | | | | | | | |
| $p=10$ | 0.9329 | 0.9396 | 0.9391 | 0.9411 | 0.9410 | 0.9464 | 0.9427 |
| $p=20$ | 0.9373 | 0.9414 | 0.9408 | 0.9448 | 0.9411 | 0.9475 | 0.9450 |
| $p=40$ | 0.9419 | 0.9482 | 0.9487 | 0.9465 | 0.9467 | 0.9451 | 0.9495 |
| $p=80$ | 0.9448 | 0.9444 | 0.9497 | 0.9496 | 0.9476 | 0.9494 | 0.9510 |
| $p=160$ | 0.9427 | 0.9413 | 0.9454 | 0.9505 | 0.9519 | 0.9473 | 0.9490 |
| $p=320$ | 0.9454 | 0.9468 | 0.9428 | 0.9451 | 0.9515 | 0.9499 | 0.9504 |
| (b) Empirical probability of accepting the identity covariance hypothesis when $\mathbf{\Sigma}_\theta = \operatorname{diag}(10, 1, \ldots, 1)$ | | | | | | | |
| $p=10$ | 0.0253 | 0.0003 | – | – | – | – | – |
| $p=20$ | 0.0531 | 0.0029 | – | – | – | – | – |
| $p=40$ | 0.1218 | 0.0093 | – | – | – | – | – |
| $p=80$ | 0.2458 | 0.0432 | 0.0080 | – | – | – | – |
| $p=160$ | 0.4263 | 0.1466 | 0.0002 | – | – | – | – |
| $p=320$ | 0.6288 | 0.3683 | 0.0858 | 0.0012 | – | – | – |

TABLE 10

*The identity covariance hypothesis is rejected when the Ledoit–Wolf test statistic in (7.6) exceeds the 5% significance level for the $\chi^2$ distribution with $p(p+1)/2$ degrees of freedom*

| | $n = 10$ | $n = 20$ | $n = 40$ | $n = 80$ | $n = 160$ | $n = 320$ | $n = 640$ |
|---|---|---|---|---|---|---|---|
| (a) Empirical probability of accepting the hypothesis $\boldsymbol{\Sigma} = \mathbf{I}$ using the Ledoit–Wolf test | | | | | | | |
| $p = 10$ | 0.9483 | 0.9438 | 0.9520 | 0.9493 | 0.9510 | 0.9553 | 0.9465 |
| $p = 20$ | 0.9498 | 0.9473 | 0.9510 | 0.9513 | 0.9498 | 0.9495 | 0.9423 |
| $p = 40$ | 0.9428 | 0.9545 | 0.9468 | 0.9448 | 0.9488 | 0.9460 | 0.9478 |
| $p = 80$ | 0.9413 | 0.9490 | 0.9513 | 0.9540 | 0.9480 | 0.9500 | 0.9460 |
| $p = 160$ | 0.9438 | 0.9495 | 0.9475 | 0.9520 | 0.9508 | 0.9543 | 0.9448 |
| $p = 320$ | 0.9445 | 0.9475 | 0.9493 | 0.9490 | 0.9485 | 0.9468 | 0.9453 |
| (b) Empirical probability of accepting the identity covariance hypothesis when $\boldsymbol{\Sigma} = \mathrm{diag}(10, 1, \ldots, 1)$ using the Ledoit–Wolf test | | | | | | | |
| $p = 10$ | 0.0345 | 0.0008 | – | – | – | – | – |
| $p = 20$ | 0.0635 | 0.0028 | – | – | – | – | – |
| $p = 40$ | 0.1283 | 0.0130 | – | – | – | – | – |
| $p = 80$ | 0.2685 | 0.0450 | 0.0008 | – | – | – | – |
| $p = 160$ | 0.4653 | 0.1575 | 0.0070 | – | – | – | – |
| $p = 320$ | 0.6533 | 0.3700 | 0.0773 | 0.0010 | – | – | – |

$p/n \to c \in (0, \infty)$,

$$(7.7) \qquad l_j \to \begin{cases} \lambda_j \left(1 + \dfrac{\lambda c}{\lambda_j - \lambda}\right), & \text{if } \lambda_j > \lambda(1 + \sqrt{c}), \\ \lambda(1 + \sqrt{c})^2, & \text{if } \lambda_j \le \lambda(1 + \sqrt{c}), \end{cases}$$

*where $j = 1, \ldots, k$ and the convergence is almost surely.*

PROOF. See Baik and Silverstein (2006), Paul (2005), Baik, Ben Arous and Péché (2005). □

Since the inference methodologies we propose in this paper exploit the distributional properties of traces of powers of the sample covariance matrix, Proposition 7.1 pinpoints the fundamental inability of the sphericity test proposed to reject the hypothesis $\boldsymbol{\Sigma} = \mathbf{I}$ whenever (for large $p, n$)

$$\lambda_i \le 1 + \sqrt{\frac{p}{n}}.$$

For the example considered, $\lambda_1 = 10$, so that the above condition is met whenever $p/n > c_t = 81$. For $p/n$ on the order of $c_t$, the resulting inability to correctly reject the identity covariance hypothesis can be attributed to this phenomenon and the fluctuations of the largest eigenvalue.

Canonically speaking, eigen-inference methodologies which rely on traces of powers of the sample covariance matrix will be unable to differentiate between closely spaced population eigenvalues in high-dimensional, sample sized starved settings. This impacts the quality of the inference in a fundamental manner that is difficult to overcome. At the same time, however, the results in Baik and Silverstein (2006), Paul (2005), Baik, Ben Arous and Péché (2005) suggest that if the practitioner has reason to believe that the population eigenvalues can be split into several clusters about $a_i \pm a_i\sqrt{p/n}$, then the use of the model in (1.2) with a block subspace structure, where the individual blocks of sizes $p_1, \ldots, p_k$ are comparable to $p$, is justified. In such situations, the benefit of the proposed eigen-methodologies will be most apparent and might motivate experimental design that ensures that this condition is met.

**8. Extensions and lingering issues.** In the development of the estimation procedures in this article, we ignored the correction term for the mean that appears in the real covariance matrix case (see Proposition 3.1). This was because Bai and Silverstein expressed it as a contour integral which appeared challenging to compute [see (1.6) in Bai and Silverstein (2004)]. It is desirable to include this extra term in the estimation procedure if it can be computed efficiently using symbolic techniques. The recent work of Anderson and Zeitouni (2006), despite its ambiguous title, represents a breakthrough on this and other fronts.

Anderson and Zeitouni encode the correction term in the coefficients of a power series that can be directly computed from the limiting moment series of the sample covariance matrix [see Theorem 3.4 in Anderson and Zeitouni (2006)]. Furthermore, they have expanded the range of the theory for the fluctuations of traces of powers of large Wishart-like sample covariance matrices, in the real sample covariance matrix case, to the situation when the entries are composed from a broad class of admissible non-Gaussian distributions. In such a scenario, the correction term takes into account the fourth moment of the distribution [see (5) and Theorems 3.3 and 3.4 in Anderson and Zeitouni (2006)]. This latter development might be of use in some practical settings where the non-Gaussianity is well characterized. We have yet to translate their results into a computational recipe for determining the correction term though we intend to do so at a later date. We plan to make a software implementation based on the principles outlined in this paper available for download. The numerical results presented show the consistency of the proposed estimators; it would be of interest to establish this analytically and identify conditions in the real covariance matrix case, where ignoring the correction term in the mean can severely degrade the quality of estimation. The issue of how a local test that exploits global information, of the sort proposed by El Karoui (2007), compares to the global test developed in this article in terms of hypothesis discriminatory power is an unresolved question of great interest. A more systematic investigation is needed of the efficacy of various model order selection techniques for the problem considered.

**Acknowledgments.** We are grateful to Debashis Paul for his feedback on an early version of this paper. His comments improved the clarity of our exposition. The initial investigation comparing the power of local versus global tests was motivated by Patrick Perry's penetrating questions; only some of the issues he raised have been addressed in this paper while others remain outstanding.

It has been a privilege to interact with Debashis Paul and Jack Silverstein on these and related matters during the authors' involvement in the SAMSI Program on High Dimensional Inference and Random Matrices. We thank Iain Johnstone for the opportunity to participate in this program and applaud James Berger, Nicole Scott, Margaret Polinkovsky, Sue McDonald and Donna Pike for creating the wonderful atmosphere at SAMSI.

We thank Boaz Nadler for detailed and helpful comments on our manuscript. The first author thanks Arthur Baggeroer for encouragement and feedback from the earliest stages of this work. We thank the referees for their excellent feedback and suggestions. In particular we appreciate a referee's insistence that we investigate the model order estimation problem further. This motivated us to streamline the numerical code so that we could run the numerical simulations reported in Section 6.2 more seamlessly.

N. R. Rao
Department of Mathematics and EECS
Massachusetts Institute of Technology
Cambridge, Massachusetts 02142
USA
E-mail: raj@mit.edu

J. Mingo
R. Speicher
Department of Mathematics
  and Statistics
Queen's University
Kingston, Ontario
Canada K7L 3N6
E-mail: mingo@mast.queensu.ca
         speicher@mast.queensu.ca

A. Edelman
Department of Mathematics
Massachusetts Institute of Technology
Cambridge, Massachusetts 02139
USA
E-mail: edelman@mit.edu